\documentclass[a4paper,10pt]{article}
\usepackage{latexsym}
\usepackage{amsmath}
\usepackage{theorem}
\usepackage{color}
\usepackage{epsfig}

\newtheorem{theorem}{Theorem}
\newtheorem{corollary}[theorem]{Corollary}
\newtheorem{definition}{Definition}
\newtheorem{lemma}[theorem]{Lemma}\newtheorem{proposition}[theorem]{Proposition} 
\newtheorem{remark}{Remark}

\newcommand{\NN}{{\rm\bf N}}

\newcommand{\ZZ}{{\rm\bf Z}}
\newcommand{\RR}{{\rm\bf R}}

\newcommand{\So}{{\mathbf {SO}}(2)}

\newcommand{\vv}{{\rm\bf v}}

\newcommand{\dpt}{\displaystyle}
\def\Qed{\hfill\raisebox{.6ex}{\framebox[2.5mm]{}} \medbreak}

\newcommand{\ord}[1]{$^{\rm \bf \underline {#1}}$}

\newcounter{lixo}
\newcounter{contai}

\title{On Takens' Last Problem: \\
tangencies and time averages near heteroclinic networks }

\date{ \today}

\author{
Isabel S. Labouriau
\quad Alexandre A. P. Rodrigues\\
Centro de Matem\'atica
da Universidade do Porto
\thanks{CMUP   (UID/MAT/00144/2013) is supported by the European Regional Development Fund through the programme COMPETE and by the Portuguese Government through the Funda\c{c}\~ao para a Ci\^encia e a Tecnologia (FCT)  under the partnership agreement PT2020. A.A.P. Rodrigues was supported by the grant
SFRH/BPD/84709/2012 of FCT. Part of this work has been written during AR stay in Nizhny Novgorod University partially supported by the grant RNF 14-41-00044.}\\
 and
Faculdade de
Ci\^encias, Universidade do Porto \\
Rua do Campo Alegre,
687, 4169-007 Porto, Portugal \\
 islabour@fc.up.pt \quad alexandre.rodrigues@fc.up.pt }
\begin{document}

\maketitle

\bigbreak

\bigbreak
\textbf{2010 --- AMS Subject Classifications} 

  {Primary: 34C28; 
    Secondary: 34C37, 37C29, 37D05, 37G35}

\bigbreak

\maketitle
\begin{abstract}
We obtain a structurally stable family of smooth ordinary differential equations exhibiting heteroclinic tangencies for a dense subset of parameters.
We use this to find vector fields $C^2$-close to an element of the family exhibiting a tangency, for which the set of solutions with historic behaviour contains an open set. 
This provides  an affirmative answer to  Taken's Last  Problem  \emph{(F. Takens (2008) Nonlinearity, 21(3) T33--T36).}
A limited solution with historic behaviour is one for which the time averages do not converge as time goes to infinity. Takens' problem asks for dynamical systems where historic behaviour occurs persistently for initial conditions in a set with  positive Lebesgue measure.

The family appears in the unfolding of a degenerate differential equation whose flow has an  asymptotically stable heteroclinic cycle involving  two-dimensional connections of non-trivial periodic solutions.
We  show that  the degenerate problem also has  historic behaviour, since  for an open set of initial conditions starting near the cycle, the time averages  approach the boundary of a polygon whose vertices depend on the  centres of gravity of the periodic solutions and their Floquet multipliers. 

We illustrate our results with an explicit example where historic behaviour arises $C^2$-close of a $\textbf{SO(2)}$-equivariant vector field. 
\end{abstract}


\section{Introduction}

Chaotic dynamics makes it difficult to give a geometric description of an attractor in many situations,
when probabilistic and ergodic analysis becomes relevant.
In a long record of a chaotic signal generated by a deterministic time evolution, for suitable initial conditions the expected time average exists --- see  \cite{Ruelle,Sigmund}. 
 However, there are cases where the time averages do not converge no matter how long we wait. 
This \emph{historic behaviour}
is associated with intermittent dynamics, which  happens typically near  heteroclinic networks.

The aim of this article is to explore the persistence of this behaviour for a deterministic class of systems involving  robust heteroclinic cycles, leading to an answer to Taken's Last   Problem \cite{T}. 
More precisely, we study non-hyperbolic heteroclinic attractors such that the time averages of all solutions within their basin of attraction do not converge, and for which this holds persistently.

This is done by first studying a one-parameter family of vector fields having periodic solutions connected in a robust cycle.
We show that under generic conditions there are parameter values for which the invariant manifolds of a pair of periodic solutions have a heteroclinic tangency.
This implies the Newhouse property of existence of infinitely many sinks. 
Results by Kiriki and Soma \cite{KS} may then be used to provide an affirmative answer to the problem proposed by Takens in \cite{T}.

\subsection{Takens' last  problem}
Let $M$ be a compact three-dimensional manifold  without boundary and consider a vector field $f: M \rightarrow TM$ defining
 a differential equation
 \begin{equation}
    \label{general}
\dot{x}=f(x), \qquad x(0)=x_0\in M
 \end{equation}
and denote by $\phi(t,x_0)$, with $t \in \RR$, the associated flow with initial condition $x_0 \in M$. 
The following terminology has been introduced by Ruelle \cite{Ruelle} (see also Sigmund \cite{Sigmund}).

\begin{definition}\label{historicDef}
We say that the solution $\phi(t,x_0)$, $x_0 \in M$, of \eqref{general} has \emph{historic behaviour} if 
there is a continuous function $H:M\rightarrow \RR $ such that
 the time average
\begin{equation}
\label{historic1}
\frac{1}{T}\int_{0}^{T} H(\phi(t,x_0)) dt
\end{equation}
fails to converge.
\end{definition}
A solution $\phi(t,x_0)$, $x_0 \in M$ with  historic behaviour  retains informations about its past.
This happens, in particular, if there are at least two different sequences of times, say $(T_i)_{i \in \NN}$ and $(S_j)_{j \in \NN}$, such that the following limits exist and are different:
$$
\lim_{i \rightarrow +\infty}\frac{1}{T_i}\int_{0}^{T_i} H(\phi(t,x_0)) dt 
\quad \neq \quad 
\lim_{j \rightarrow +\infty}\frac{1}{S_j}\int_{0}^{S_j} H(\phi(t,x_0)) dt.
$$

The consideration of the  limit behaviour  of time averages with respect to a given measure has been studied
since Sinai \cite{Sinai}, Ruelle \cite{Ruelle76} and Bowen \cite{Bowen75}.
Usually, historic behaviour is seen as an anomaly.
Whether there is a justification for this belief is the content of Takens' Last Problem \cite{KS, Takens94, T}: \emph{are there persistent classes of smooth dynamical systems such that the set of initial conditions which give rise to orbits with historic behaviour has positive Lebesgue measure?} In ergodic terms, this problem is equivalent to finding a persistent class of systems admitting no physical measures \cite{Hofbauer, Ruelle}, since roughly speaking, these measures are those that give probabilistic information on the observable asymptotic behaviour of trajectories.

 The class may become persistent if one considers differential equations in manifolds with boundary as in population dynamics \cite{Hofbauer, HSig}. The same happens for equivariant or reversible differential equations \cite{Guckenheimer e Holmes 1}. The question remained open for systems without such properties
until, recently, Kiriki and Soma \cite{KS}  proved that any Newhouse open set in the $C^r$-topology,  $r\geq 2$, of two-dimensional diffeomorphisms is contained in the closure of the set of diffeomorphisms which have non-trivial wandering domains whose forward orbits have historic behaviour.  As far as we know, the original problem, stated   for flows, has remained open until now.

\subsection{Non-generic historic behaviour}

In this section, we present some non-generic examples that, however,
occur generically in families of discrete dynamical systems depending on a small number of parameters.
The first example has been given in Hofbauer and Keller \cite{HK}, where it has been shown that the logistic family contains elements for which almost all orbits have historic behaviour. This example has codimension one in the space of $C^3$ endomorphisms of the interval; the $C^3$ regularity is due to the use of the Schwarzian derivative operator.

The second example is due to Bowen, who  described a  codimension two
system of differential equations on the plane whose flow has a heteroclinic cycle
 consisting of a pair of saddle-equilibria connected by two trajectories. As referred by Takens \cite{Takens94, T}, apparently Bowen never published this result. 
 We give an explicit example in \ref{subsecBowen} below.
 The eigenvalues of the derivative of the vector field at the two saddles are such that the cycle attracts solutions that start inside it. In this case, each solution in the domain has historic behaviour.
 In ergodic terms,  it is an example without SRB measures.
  Breaking the cycle by a small perturbation, the  equation loses this property. This type of dynamics may become persistent for dynamical systems in manifolds with boundary or in the presence of symmetry. 
 We use Bowen's example here as a first step in the construction of a generic example.
  Other examples of high codimension with heteroclinic attractors where Lebesgue almost all trajectories fail to converge have been given by Gaunersdorfer \cite{Gaunersdorfer} and Sigmund \cite{Sigmund}.

Ergodicity implies the convergence of time averages along almost all trajectories for all continuous observables \cite{KA}. 
For non-ergodic systems, time averages may not exist for almost all trajectories. 
 In Karabacak and Ashwin \cite[Th 4.2]{KA}, the authors characterise conditions on the observables that imply convergent time averages for almost all trajectories. This convergence is determined by the behaviour of the observable on the statistical attractors (subsets where trajectories spend almost all time). Details in \cite[\S 4]{KA}.

\subsection{General examples}
The paradigmatic example with persistent historic behaviour has been suggested by Colli and Vargas in \cite{CV}, 
in which the authors presented a simple non-hyperbolic model with a wandering domain characterised by the existence  of a two-dimensional diffeomorphism with a Smale horseshoe whose stable and unstable manifolds have persistent tangencies under arbitrarily small $C^2$ perturbations. The authors of \cite{CV} suggest that this would entail the existence of non-wandering domains with historic behaviour, in a robust way. This example has been carefully described in \cite[\S 2.1]{KS}.

For diffeomorphisms, an answer has been given by Kiriki and Soma \cite{KS}, where the authors used ideas suggested in \cite{CV} to find  a nontrivial non-wandering domain (the interior of a specific rectangle) where the diffeomorphism is contracting. In a robust way, they obtain an open set of initial conditions for which the time averages do not converge.
Basically, the authors linked two subjects: homoclinic tangencies studied by Newhouse, Palis and Takens and non-empty non-wandering domains exhibiting historic behaviour.   
An overview of the proof has been given in \S 2 of \cite{KS}. 
We refer those that are unfamiliar with Newhouse regions to the book \cite{PT}.

\subsection{The results}
The goal of this article is twofold. First, we  extend the results by Takens \cite{Takens94} and by Gaunersdorfer \cite{Gaunersdorfer}  to heteroclinic cycles involving periodic solutions with real Floquet multipliers. The first main result is Theorem~\ref{Main1}, with precise hypotheses   given in Section \ref{Hypotheses}:
\begin{description}
\item[1\ord{st} result:] 
Consider an ordinary differential equation in $\RR^3$ having  an attracting heteroclinic cycle involving periodic solutions with  two-dimensional heteroclinic connections. 
Any neighbourhood of this cycle contains an open set of initial conditions, for which the  time averages of the corresponding solutions  accumulate on the boundary of a polygon, and thus,  fail to converge.
The open set is contained in the basin of attraction of the cycle and the observable is the projection on a component.
\end{description}
 This situation has  high codimension because each heteroclinic connection raises the codimension by one, but this class of systems  is persistent in equivariant differential equations. The presence of  symmetry creates
flow-invariant fixed-point subspaces in which heteroclinic connections lie --- see for example the example constructed in \cite[\S 8]{Rodrigues}. 
Another example is constructed in Section~\ref{secLifting} below. The second main result, Theorem~\ref{teorema tangency}, concerns tangencies:
\begin{description}
\item[2\ord{nd} result:]  
Consider a generic one-parameter family of structurally stable differential equations in the unfolding of an equation for which the 1$^{\rm\underline{ st}}$ result holds.
Then  there is a sequence of parameter values for which there is a heteroclinic tangency of the invariant manifolds of two periodic solutions.

\end{description}

We use this result to obtain   Theorem~\ref{teoremaHistoric}:
\begin{description}
\item[3\ord{rd} result:]  
Consider a generic one-parameter family of structurally stable differential equations in the unfolding of an equation for which the 1$^{\rm\underline{ st}}$ result holds. 
Therefore, for parameter values in an open interval, there are vector fields arbitrarily $C^2$-close to an element of the family, for which there is an open set of initial conditions exhibiting historic behaviour.
\end{description}
In other words, we obtain a class, dense in a  $C^2$-open set of differential equations and elements of this class exhibit historic behaviour for an open set  of initial conditions, which may be interpreted as the condition required in Takens' Last Problem.
The idea behind the proof goes back to the works of \cite{LR3,LR2015}, combined with the recent progress on the field made by \cite{KS}. The proof consists of the followingsteps:
\begin{enumerate}
 \item  use the 3\ord{rd} result to establish the existence of intervals in the parameters corresponding to Newhouse domains;
\item\label{Item1}  in a given cross section, construct a diffeomorphim ($C^2$-close to the first return map) having historic behaviour for an open set of initial conditions; 
\item\label{Item2} 
transfer the historic behaviour from the perturbed diffeomorphism of \emph{\ref{Item1}.} to  a flow  $C^2$-close to the original one.
\end{enumerate}
Furthermore, in the spirit of the  example by  Bowen described in \cite{Takens94}, we obtain:
\begin{description}
\item[4\ord{th} result:]  We  construct  explicitly a class of systems for which historic behaviour arises $C^2$-close to the unfolding of a fully symmetric vector field, we may find an open set of initial conditions with historic behaviour.  
In contrast to the findings of Bowen and Kleptsyn \cite{Kleptsyn}, our example is robust due the hyperbolicity of the periodic solutions and the transversality of the local heteroclinic  connections. 
\end{description}

The results in this article are stated for vector fields in $\RR^3$, but they hold for vector fields in  a three-dimensional Riemannian manifold and, with some adaptation,  in higher dimensions.

\subsection{An ergodic point of view}
Concerning the first result, the outstanding fact  in the degenerate case is that the time averages diverge precisely in the same way: they approach a $k$-polygon. This is in contrast with ergodic and hyperbolic strange attractors admitting a physical measure, where almost all initial conditions lead to converging time averages, in spite of the fact that
 the observed dynamics may undergo huge variations.
 
 If a flow $\phi(t, .)$ admits an invariant probability measure $\mu$ that is absolutely continuous with respect to the Lebesgue measure and ergodic, then $\mu$ is a physical measure for $\phi(t, .)$, as a simple consequence of the  Birkhoff Ergodic Theorem. 
 In other words  if $H: M \rightarrow \RR$ is a $\mu$-integrable function, then for $\mu$-almost all points in $M$ the time average:
$$
\lim_{T \rightarrow + \infty}\frac{1}{T}\int_{0}^{T} H\circ \phi(t,x_0) dt
$$
exists and equals the space average $\int H d\mu$. In the conservative context, historic behaviour has zero Lebesgue measure.

 Physical measures need not be unique or even exist in general. 
 When they exist, it is desirable that the set of points whose asymptotic time averages are described by physical measures be of full Lebesgue measure. 
 It is unknown in how much  generality do the basins of physical measures cover   a  subset of $M$ of full Lesbegue measure.
 There are examples of systems admitting no physical measure but the only known cases are not robust, \emph{ie}, there are systems arbitrarily close (in the $C^2$ Whitney topology)  that admit physical measures. In the present  article, we exhibit a persistent class of smooth dynamical systems that does not have global physical measures.
 In the unfolding  of an equation for which the first result holds, there are no physical measures whose basins intersect the basin of attraction of an attracting heteroclinic cycle. Our example confirms that physical measures need not exist for all vector fields. Existence results are usually difficult and are known only for certain classes of systems. 

\subsection {Example without historic behaviour}
Generalised Lotka-Volterra systems has been analysed by Duarte \emph{et al} in \cite{DFO}.
Results about the convergence of time averages are known in two cases: either if there exists a unique interior equilibrium point, or in the conservative setting  (see  \cite{DFO}), when there is a heteroclinic cycle.
In the latter case, if the solution remains limited and does not converge  to the cycle, then its time averages converge to an equilibrium point. The requirement is that the heteroclinic cycle is stable but not attracting, and the limit dynamics has been extended to polymatrix replicators in \cite{Peixe}.
This is in contrast to our findings in the degenerate case, emphasising the importance of the hypothesis that the cycle is attracting in order to obtain convergence to a polygon.

\subsection{Framework of the article}

Preliminary definitions are the subject of Section~\ref{Preliminaries} and the main hypotheses are stated in Section~\ref{Hypotheses}.
We introduce the notation for the rest of the article in Section~\ref{Local} after a linearisation of the vector field around each periodic solution, whose details are given in  Appendix~\ref{appendix}.
 We use precise control of  the times of flight between cross-sections in Section~\ref{Organizing}, to show that
 for an open set of initial conditions  in a neighbourhood of asymptotically stable heteroclinic cycles involving non-trivial periodic solutions, the time averages fail to converge.
Instead, the time averages accumulate on  the boundary of a polygon, whose vertices may be computed from local information on the periodic solutions in the cycle.
The proofs of some technical lemmas containing the computations about the control of the flight time between nodes appear in Appendix \ref{appendixB}, to make for easier reading.

In Section~\ref{Tangencies}, we obtain a persistent class of smooth dynamical systems such that an open set of initial conditions corresponds to trajectories with historic behaviour.
Symmetry-breaking techniques are used to obtain a heteroclinic cycle associated to two periodic solutions and we find heteroclinic tangencies and 
Newhouse phenomena near which the result of \cite{CV, KS} may be applied.  
This is followed in Section \ref{Example} by an explicit example where historic behaviour arise in the unfolding of an $\textbf{SO(2)}$-equivariant vector field.

 \section{Preliminaries}
\label{Preliminaries}
 To make the paper self-contained and readable, we recall some definitions. 
\subsection{Heteroclinic attractors} 

Several definitions of heteroclinic cycles and networks have been given in the
literature. In this paper we consider non-trivial periodic solutions of (\ref{general}) that are hyperbolic and that have one Floquet multiplier with absolute value greater than 1 and one Floquet multiplier with absolute value less than 1. 
A connected component of $W^s(\mathcal{P})\backslash \mathcal{P}$, for  a periodic solution $\mathcal{P}$, will be called a \emph{branch} of $W^s(\mathcal{P})$, with a similar definition for a branch of $W^u(\mathcal{P})$.
Given two periodic solutions $\mathcal{P}_{a}$ and $\mathcal{P}_{b}$ of (\ref{general}), a \emph{heteroclinic connection} from $\mathcal{P}_a$ to $\mathcal{P}_b$ is a trajectory
contained in $W^u(\mathcal{P}_a)\cap W^s(\mathcal{P}_b)$,  that will be denoted $[\mathcal{P}_a\to  \mathcal{P}_b]$. 

Let $\mathcal{S=}\{\mathcal{P}_{a}:a\in \{1,\ldots,k\}\}$ be a finite ordered set of periodic solutions of saddle type of (\ref{general}).
 The notation for $\mathcal{P}_{a}$ is  cyclic, we indicate this by taking the index $a\pmod{k}$, \emph{ie} $a\in\ZZ_k = \ZZ/k\ZZ$. Suppose
 $$
 \forall a\in\ZZ_k
  \quad W^{u}(\mathcal{P}_{a})\cap W^{s}(\mathcal{P}_{a+1})\neq\emptyset .
 $$
A \emph{heteroclinic cycle} $\Gamma$ associated to $\mathcal{S}$ is the union of the saddles in 
$\mathcal{S}$ with a heteroclinic connection 
$[\mathcal{P}_a \rightarrow \mathcal{P}_{a+1}]$ 
for each $a\in\ZZ_k$.
We refer to the saddles defining the heteroclinic cycle as \emph{nodes}.  
A \emph{heteroclinic network} is a connected set that is the union of heteroclinic cycles.
When a branch of $W^{u}(\mathcal{P}_{a})$ coincides with a branch of $W^{s}(\mathcal{P}_{a+1})$, we also refer to it as a two-dimensional connection $[\mathcal{P}_a \rightarrow \mathcal{P}_{a+1}]$.

\subsection{Basin of attraction}
For a solution of (\ref{general}) passing through $x\in M$, the set of its accumulation points as $t$ goes to $+\infty$ is the $\omega$-limit set of $x$ and will be denoted by $\omega(x)$. More formally, 
$$
\omega(x)=\bigcap_{T=0}^{+\infty} \overline{\left(\bigcup_{t>T}\phi(t, x)\right)}.
$$ 
It is well known that $\omega(x)$  is closed and flow-invariant, and if $M$ is compact, then $\omega(x)$ is non-empty for every  $x\in M$.
If $\Gamma\subset M$ is a flow-invariant subset for (\ref{general}), the \emph{basin of attraction of $\Gamma$} is given by
 $$
 \mathcal{B}(\Gamma) = \{x \in M\backslash \Gamma : \mbox{all accumulation points of } \phi(t, x)\mbox{ as } t\to +\infty \mbox{ lie in } \Gamma\} .
 $$
Note that, with this definition, the set $\Gamma$ is not contained in $\mathcal{B}(\Gamma)$.

\section{The setting}
\label{Hypotheses}
\subsection{The hypotheses}
\begin{figure}
\begin{center}
\includegraphics[height=4.5cm]{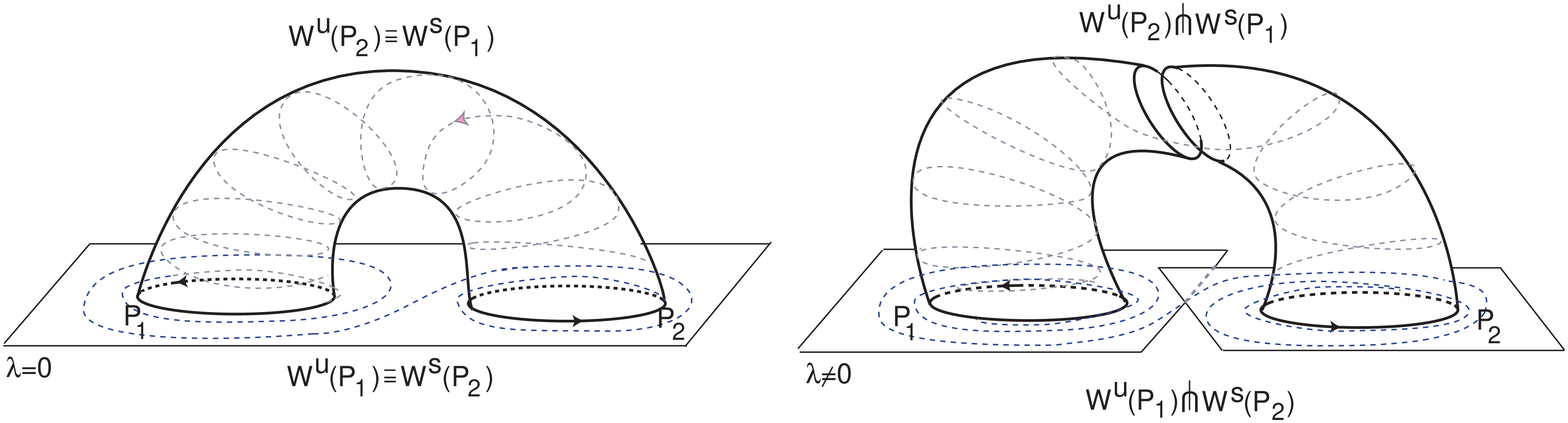}
\end{center}
\caption{\small Configuration of $\Gamma_\lambda$ for $\lambda=0$ (left) and $\lambda>0$ (right). The representation is done for $k=2$.}
\label{Configuration}
\end{figure}

Our object of study is the dynamics around a heteroclinic cycle associated to $k$ periodic solutions,
$k\in\NN$, $k>1$,  for which we give a rigorous
description here. 
Specifically, we study a one-parameter family of $C^2$-vector fields $f_\lambda$ in $\RR^3$ whose flow has the following properties (see Figure \ref{Configuration}):

\begin{enumerate}
\renewcommand{\labelenumi}{(P{\theenumi})}
\item\label{P1} 
For $\lambda\in \RR$, there are $k$ hyperbolic periodic solutions $\mathcal{P}_a$ of $\dot{x}=f_\lambda(x)$,
$a\in\ZZ_k$,
 of minimal period $\xi_a>0$.  The Floquet multipliers of $\mathcal{P}_a$ are real and given by $e^{e_a}>1$ and $e^{-c_a}<1$ where $c_a> e_a>0 $.
\item\label{P2} 
For each $a\in\ZZ_k$,
 the manifolds $W^s_{loc}(\mathcal{P}_a)$ and $W^u_{loc}(\mathcal{P}_a)$ are smooth surfaces homeomorphic to a cylinder -- see Figure \ref{local_C}.
\item\label{P3} 
For each $a\in\ZZ_k$, and for $\lambda=0$, one branch of  $W^u(\mathcal{P}_{a})$ coincides with a branch of
$W^s(\mathcal{P}_{a+1})$,  forming a heteroclinic  network, that we call $\Gamma_0$, and whose basin of attraction contains an open set.
\item\label{P4}[Transversality]
For $\lambda\neq 0$ and for each $a\in\ZZ_k$, a branch of the two-dimensional manifold 
$W^u (\mathcal{P}_{a})$ 
intersects transverselly a branch of
 $W^s(\mathcal{P}_{a+1})$ 
 at two trajectories, forming a heteroclinic  network $\Gamma_\lambda$, consisting of two heteroclinic cycles.
\setcounter{lixo}{\value{enumi}}
\end{enumerate}
For $\lambda\neq 0$, any one of the two trajectories of (P\ref{P4}) in $W^u (\mathcal{P}_{a})\cap W^s(\mathcal{P}_{a+1})$ will be denoted by $[\mathcal{P}_{a}\to \mathcal{P}_{a+1}]$. A more technical assumption (P\ref{P5}) will be made in Section~\ref{subsecSuspension} below, after we have established some notation.  For $a\in\ZZ_k$,  define the following constants: 
 \begin{equation}\label{constants}
 \delta_a=\frac{c_a}{e_a} >1, \qquad \mu_{a+1}= \frac{c_a}{e_{a+1}} \qquad \text{and} \qquad \delta=\prod_{a=1}^k \delta_a >1
 \end{equation}

 Also denote by $\overline{x}_a\in \RR^4$ the centre of gravity of $\mathcal{P}_{a}$,  given by $$\overline{x}_a=\frac{1}{\xi_a} \int_{0}^{\xi_a}\mathcal{P}_{a} (t) dt \in \RR^3.$$

 Without loss of generality we assume that the minimal period $\xi_a=1$, for all $a\in\ZZ_k$.
 It will be explicitly used  in system (\ref{ode of suspension}) below.

\subsection{The dynamics}
The dynamics of this kind of heteroclinic structures involving periodic solutions has been studied before in \cite{ACL NONLINEARITY, ALR, Melbourne, Rodrigues}, in different contexts.

Since $f_0$ satisfies  (P\ref{P1})--(P\ref{P3}) then, adapting the Krupa and Melbourne criterion \cite{KM1, KM2}, any solution starting sufficiently close to $\Gamma_0$ will approach it in positive time; in other words  $\Gamma_0$ is asymptotically stable.  
 As a trajectory approaches $\Gamma_0$,  it
 visits one periodic solution, then moves off to visit the other periodic solutions in the  network.
 After a while it returns to visit  the initial periodic solution, and the second  visit lasts longer than the first.
 The oscillatory regime of such a solution seems to switch into different nodes, at geometrically increasing times. 

For $\lambda \neq 0$,  by (P\ref{P4}), the invariant manifolds of the nodes meet transversally, and the network is no longer asymptotically stable due to the presence of suspended horseshoes in its neighbourhood.  As proved in \cite{Rodrigues}, there is an infinite number of heteroclinic and homoclinic connections between any two periodic  solutions and the dynamics near the heteroclinic network is very complex. The route to  chaos corresponds to an interaction of robust switching with chaotic cycling. The emergence of chaotic cycling does not depend on the magnitude of the multipliers of the periodic solutions. It depends only on the geometry of the flow near the cycle.   
 \medbreak
  In Table \ref{notation}, we summarise some information about the type of heteroclinic structure of  $\Gamma_\lambda$ and the type of dynamics nearby. 
\begin{table}[ht]
\begin{center}
\begin{tabular}{llll}
$\lambda$ & 
Structure of $V_{\Gamma_\lambda}$ & Dynamics near $\Gamma_\lambda$  & References \\ \hline
zero& 
torus of genus $k$ & Attractor & \cite{Melbourne, Rodrigues}\\ \hline
non-zero & 
torus of genus $> k$  &Chaos (Switching and Cycling) & \cite{ACL NONLINEARITY, ALR, Rodrigues}\\ \hline
\end{tabular}
\end{center}
\caption{ \small Heteroclinic structure of  $\Gamma_\lambda$, for $\lambda=0$ and $\lambda\neq0$.} 
\label{notation}
\end{table}

\section{Local and global dynamics near the network}
\label{Local}
Given a heteroclinic network of periodic solutions $\Gamma_\lambda$ with nodes $\mathcal{P}_{a}$, $a\in\ZZ_k$, let $V_{\Gamma_\lambda}$ be a compact neighbourhood of $\Gamma_\lambda$ and let $V_a$
be pairwise disjoint compact neighbourhoods of the nodes  $\mathcal{P}_{a}$, such that each boundary  $\partial V_a$ is a finite union of smooth manifolds 
with boundary, that are transverse to the vector field everywhere, except at their boundary.  
Each  $V_a$ is called an \emph{isolating block} for  $\mathcal{P}_{a}$ and, topologically, it consists of a hollow cylinder.  Topologically, $V_{\Gamma_0}$ may be seen as a   solid torus with  genus $k$ (see Table~\ref{notation}).

\subsection{Suspension and local coordinates}\label{subsecSuspension}
For $a\in\ZZ_k$, let $\Sigma_a$ be a cross section transverse to the flow at $p_a \in \mathcal{P}_{a}$. Since $\mathcal{P}_{a}$ is hyperbolic, there is a neighbourhood $V^*_a$ of $p_a$ in $\Sigma_a$ where the first return map to $\Sigma_a$, denoted by $\pi_a$, is $C^1$ conjugate to its linear part. 
Moreover, for each $r\ge 2$ there is an open and dense subset of $\RR^2$ such that, if the eigenvalues $(c_a,e_a)$ lie in this set, then the conjugacy is of class $C^r$ --- see \cite{Takens71} and Appendix~\ref{appendix}.
The eigenvalues of $d\pi_a$ are $e^{e_a}$ and $e^{-c_a}$. Suspending the linear map gives rise, in cylindrical coordinates $(\rho, \theta, z)$ around $\mathcal{P}_{a}$, to the system of differential equations:
\begin{equation}
\label{ode of suspension}
\left\{ 
\begin{array}{l}
\dot{\rho}=-c_{a}(\rho -1) \\ 
\dot{\theta}=1 \\ 
\dot{z}=e_{a}z
\end{array}
\right.
\end{equation}
which is $C^2$-conjugate, after reparametrising the time variable, to the original flow near $\mathcal{P}_{a}$. In these coordinates, the periodic solution $\mathcal{P}_{a}$ is the circle defined by $\rho=1$ and $z=0$, 
its local stable manifold, $W^s_{loc}(\mathcal{P}_{a})$, is  the plane defined by $z=0$ and $W^u_{loc}(\mathcal{P}_{a})$ is the surface defined by $\rho=1$ as  in  Figure~\ref{local_C}.

We will work with a hollow three-dimensional cylindrical neighbourhood $V_a(\varepsilon)$ of 
$\mathcal{P}_{a}$ contained in the suspension of $V^*_a$  given by:
$$
V_a(\varepsilon)=\left\{ (\rho,\theta,z):\quad 1-\varepsilon\le\rho\le 1+\varepsilon,
\quad -\varepsilon\le z\le \varepsilon\quad \text{and}\quad 
\theta\in\RR\pmod{2\pi}
\right\}\  .
$$
When there is no ambiguity, we write $V_a$ instead of $V_a(\varepsilon)$.
Its boundary is a disjoint union 
$$
\partial V_{a}= In(\mathcal{P}_{a}) \cup Out(\mathcal{P}_{a}) \cup \Omega(\mathcal{P}_{a})
$$
such that :
\begin{itemize}
\item 
$In(\mathcal{P}_{a})$ is the union of the  walls, defined by $\rho=1\pm\varepsilon$, of the cylinder,
 locally separated by $W^u(\mathcal{P}_{a})$.
Trajectories starting at  $In(\mathcal{P}_{a})$ go inside the cylinder $V_a$ in small positive
time.
\item 
$Out(\mathcal{P}_{a})$ is the union of  two anuli, the top and the bottom of the cylinder, defined by $z=\pm\varepsilon$, locally separated by $W^s(\mathcal{P}_{a})$.
Trajectories starting at $Out(\mathcal{P}_{a})$ go inside the cylinder $V_a$ in small
negative time.
\item 
The vector field is transverse to $\partial V_{a}$ at all points except possibly at the
four circles:
$\Omega(\mathcal{P}_{a})=\overline{In(\mathcal{P}_{a})}\cap \overline{Out(\mathcal{P}_{a})}$.
\end{itemize}

The two cylinder walls, $In(\mathcal{P}_{a})$ are
parametrised by the covering maps:
$$
(\theta,z)\mapsto(1\pm\varepsilon,\theta,z)=(\rho,\theta,z),
$$
where $\theta\in\textbf{R}\pmod{2\pi}$, $|z|<\varepsilon$. 
In these coordinates, $In(\mathcal{P}_{a})\cap W^s(\mathcal{P}_{a})$ is the union of the two circles 
$z=0$. 
The two anuli $Out(\mathcal{P}_{a})$ are parametrised by the coverings:
$$
(\varphi,r) \mapsto ( r,\varphi, \pm \varepsilon)=(\rho,\theta,z),
$$
for $1-\varepsilon<r<1+\varepsilon$ and $\varphi \in \RR\pmod{2\pi}$ and
where $Out(\mathcal{P}_{a})\cap W^u(\mathcal{P}_{a})$ is the union of the two circles  $r=1$.
In these coordinates $\Omega(\mathcal{P}_{a}) =\overline{In(\mathcal{P}_{a})}\cap \overline{Out(\mathcal{P}_{a})}$ is the union of the four circles defined by 
$\rho=1\pm \varepsilon$ and $ z=\pm \varepsilon$.

The portion of the unstable manifold of $\mathcal{P}_{a}$  that goes from $\mathcal{P}_{a}$ to $In(\mathcal{P}_{a+1})$ without intersecting $V_{a+1}$ will be  denoted $W^u_{loc}(\mathcal{P}_{a})$. Similarly, $W^s_{loc}(\mathcal{P}_{a})$ will denote the portion of the stable manifold of $\mathcal{P}_{a}$ that is outside $V_{a-1}$ and goes directly from $Out(\mathcal{P}_{a-1})$ to $\mathcal{P}_{a}$.
With this notation, we formulate the following technical condition:
\begin{figure}
\begin{center}
\includegraphics[height=8cm]{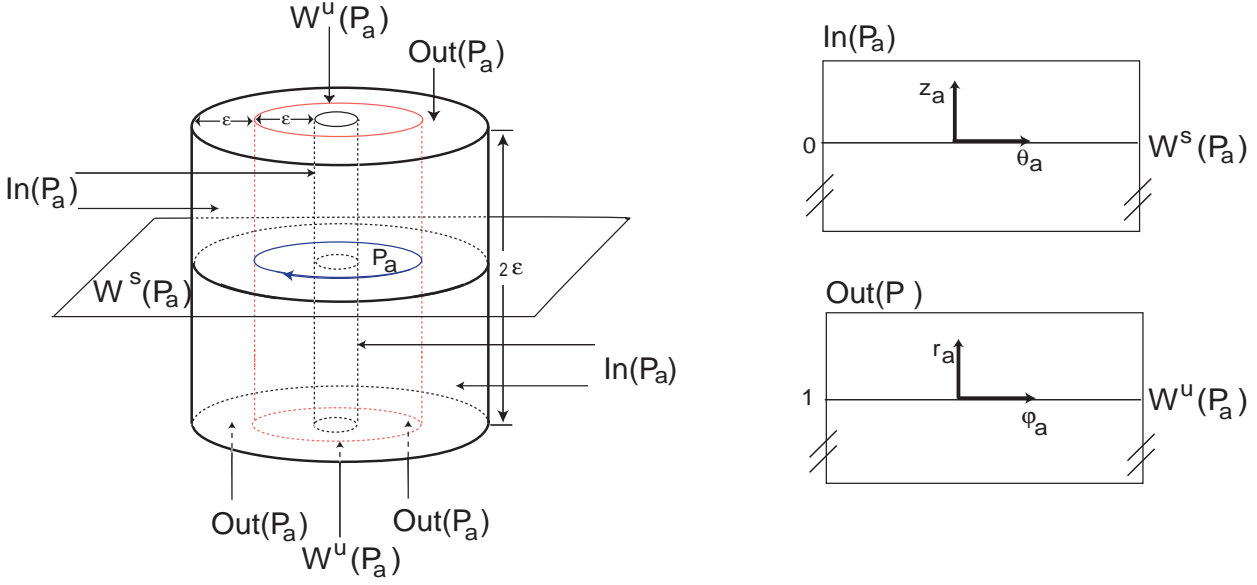}
\end{center}
\caption{\small Local coordinates on the boundary of the neighbourhood $V_a$ of a periodic solution $\mathcal{P}_a$ where $a \in \ZZ_k$. Double bars mean that the sides are identified.}
\label{local_C}
\end{figure}

 \begin{enumerate}
\renewcommand{\labelenumi}{(P{\theenumi})}
\setcounter{enumi}{\value{lixo}}
\item\label{P5} For $a\in\ZZ_k$,
and $\lambda\neq 0$ close to zero, the manifolds $W^u_{loc}(\mathcal{P}_{a})$ intersect the cylinders $In(\mathcal{P}_{a+1})$  on a closed curve. 
Similarly, $W^s_{loc}(\mathcal{P}_{a})$ intersects the  annulus $Out(\mathcal{P}_{a-1})$  on a closed curve.
\setcounter{lixo}{\value{enumi}}
\end{enumerate}

The previous hypothesis complements (P\ref{P4}) and corresponds to the expected
unfolding from the coincidence of the manifolds $W^s(\mathcal{P}_{a+1})$ and $W^u(\mathcal{P}_{a})$ at $f_0$, see Chilingworth \cite{Chilingworth}.
Note that (P\ref{P4}) and  (P\ref{P5}) are satisfied in an open subset of the set of unfoldings $f_\lambda$ of $f_0$ satisfying (P\ref{P1})--(P\ref{P3}). 

  In order to distinguish the local coordinates near the periodic solutions, we sometimes add the index $a$ with $a\in\ZZ_k$.

\subsection{Local map near the periodic solutions}\label{sublocal}
For each $a \in \ZZ_k$, we may solve  (\ref{ode of suspension}) explicitly, then we compute the flight time from  $In(\mathcal{P}_{a})$ to $Out(\mathcal{P}_{a})$ by solving the equation $z(t)=\varepsilon$ for the trajectory whose initial condition is $(\theta_a, z_a) \in In(\mathcal{P}_{a})\backslash W^s(\mathcal{P}_{a})$, with $z_a>0$.
We find that this trajectory  arrives at $ Out(\mathcal{P}_{a})$ at a time $\tau_a: In(\mathcal{P}_{a})\backslash W^s(\mathcal{P}_{a}) \rightarrow \RR_0^+$ given by:
 \begin{equation}
 \label{Time of Flight}
   \tau_a(\theta_a, z_a)=\frac{1}{e_a}\ln \left(\frac{\varepsilon}{z_a}\right).
 \end{equation}
 Replacing this time in the other coordinates of the solution, yields:

  \begin{equation}
\label{local map}
\Phi _{a}(\theta_a,z_a)=
\left(\theta_a-\frac{1}{e_a}\ln\left(\frac{z_a}{\varepsilon}\right),
1\pm \varepsilon \left(\frac{z_a}{\varepsilon}\right)^{\delta_a}\right)= (\varphi_a,r_a) 
\qquad\mbox{where} \quad\delta_a=\frac{c_{a}}{e_{a}}>1 .
\end{equation}
The signs $\pm$ depend on the component of $In(\mathcal{P}_{a})$ we started at,
 $+$ for trajectories starting with $r_a>1$ and $-$ for  $r_a<1$. 
 We will discuss  the case $r_a>1$, $z_a>0$, the behaviour on the other components is analogous.

\subsection{Flight times for $\lambda=0$}
\label{secTransition0}
 
\begin{figure}
\begin{center}
\includegraphics[height=7.5cm]{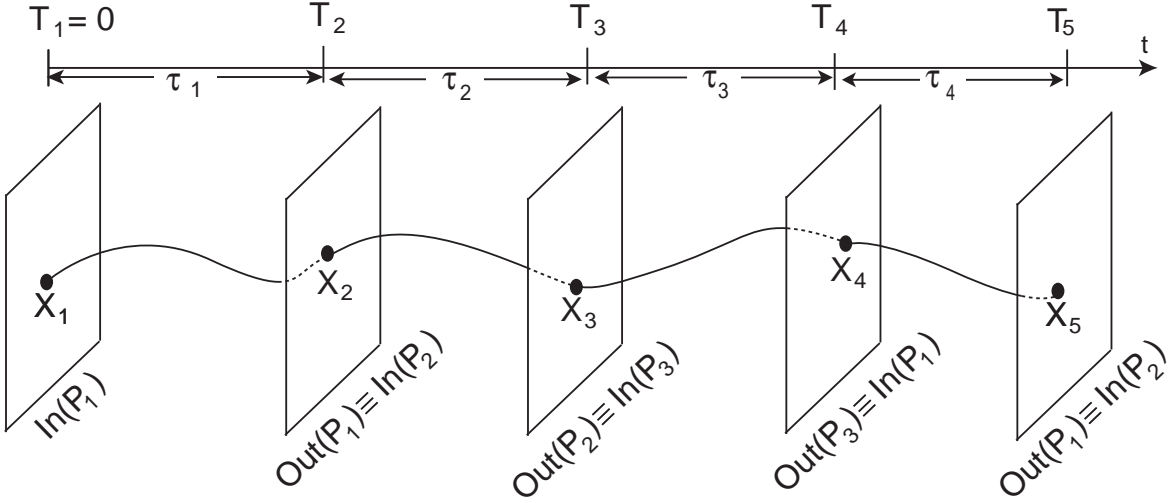}
\end{center}
\caption{\small For $X\in \mathcal{B}(\Gamma_0)$, the solution $\phi(t,X)$  remains in ${V}_a$ for a time interval of length $\tau_{a}(X)$, then spends $\tau_{a+1}(X)$ units of time near $\mathcal{P}_{a+1}$, and ,
after $n$ full turns, stays again in ${V}_a$ for $\tau_{a+nk}(X)$ units of time, and so on. The representation is done for $k=3$.}
\label{times_notation}
\end{figure}

 Here we introduce some terminology that will be used in Section~\ref{Organizing};  see Figure \ref{times_notation}.
 For $X\in \mathcal{B}(\Gamma_0)$, let $T_1(X)$ be the smallest $t\ge 0$ such that  $\phi(t,X)\in In(\mathcal{P}_{1})$.
 For $j\in \NN$, $j>1$, we define $T_j(X)$ inductively as the smallest $t>T_{j-1}(X)$ such that $\phi(t,X)\in In(\mathcal{P}_{\langle j\rangle})$, where 
$$ 
\left\langle j\right\rangle= j- \left[\frac{j}{k}\right]k 
$$ 
 is the remainder in the integer division by $k$ and $[x]$ is the greatest integer less than or equal to $x$. Recall that the index $a$ in $\mathcal{P}_{a}$ lies in $\ZZ_a$, so that $\mathcal{P}_{0}$ and $\mathcal{P}_{k}$ represent the same periodic solution.
 
 In order to simplify the computations, we may assume that the transition from $Out(\mathcal{P}_a)$ to $In(\mathcal{P}_{a+1})$ is instantaneous.
This is reasonable because, as $t\to\infty$, the time of flight inside each $V_a$ tends to infinity, whereas the time of flight from $Out(\mathcal{P}_a)$ to $In(\mathcal{P}_{a+1})$ remains limited.
In the proof of Proposition~\ref{density} below, we will see that this assumption does not affect the validity of our results. 
With this assumption, the time of flight $\tau_{a+nk}(X)$ inside $V_a$ at the $n$-th pass of the trajectory through $V_a$ will be 
$$
\tau_{a+nk}(X)=T_{a+1+nk}(X)-T_{a+nk}(X),
$$
thus extending the notation $\tau_a$ introduced in \ref{sublocal} above   to $X\in \mathcal{B}(\Gamma_0)$ and any index $a+nk\in\NN$.

\medbreak
For each $a \in \ZZ_k$, and for $\lambda=0$, we define the transition map
$ \Psi_{a }^0:Out(\mathcal{P}_{a})\rightarrow In(\mathcal{P}_{a+1})$ \begin{equation}
\label{Psi_def}
\Psi_a^0(\varphi_a, r_a) = (\varphi_a, r_a-1)=({\theta}_{a+1}, {z}_{a+1}).
\end{equation}
The transition maps for $\lambda\ne 0$ will being discussed in Section~\ref{secTransition}.
\medbreak

\section{The $k$-polygon at the organising centre}
\label{Organizing}	
 Let $f_0$ be a vector field in $\RR^3$ satisfying (P\ref{P1})--(P\ref{P3}). All the results of this section assume $\lambda=0$. 
  Suppose, from now on, that $\phi(t, X) $ is a solution $\dot{x}=f_0(x)$ with initial condition $X=\phi(0,X)$  in $\mathcal{B}(\Gamma_0)$, the basin of attraction of  $\Gamma_0$.

\subsection{The statistical limit set of $f_0$}

The statistical limit set $\Lambda_{stat}(f_0)$ associated to the basin of attraction of $\Gamma_0$ is the smallest closed subset where Lebesgue almost all trajectories spend almost all time. More formally, following Ilyashenko \cite{Ilya1} and Karabacak and Ashwin \cite{KA}, we define: 

\begin{definition} For an open set $U\subset \RR^3$ and a solution $\phi(t,x)$ of (\ref{general}) with $x\in \RR^3$:
\begin{enumerate}
\item the frequency of the solution being in $U$ is the  ratio: 
$$
\rho_{f}(x, U, T)= \frac{Leb\{t \in  [0,T]: \phi(t,x) \in U\}}{T}.
$$
where $Leb$ denotes the Lebesgue measure in $\RR$.
\item the statistical limit set, denoted by $\Lambda_{stat}({f})$, is the smallest closed subset of $\RR^3$ for which any open neighbourhood of $U$ of $\Lambda_{stat}$ satisfies the equality:
$$
\lim_{t\rightarrow +\infty} \rho_{f}(x, U, t)=1, \qquad \text{for almost all  } x \in \RR^3.
$$ 
\end{enumerate}
\end{definition}
Since the transitions between the saddles of $\Gamma_0$ are very fast compared with the times of sojourn near the periodic solutions $\mathcal{P}_a$, $a\in\ZZ_k$  (see \ref{secTransition0}) we may conclude that: 
\begin{proposition}
\label{density}
Let $f_0$ be a vector field  in $\RR^3$  satisfying (P\ref{P1})--(P\ref{P3}). Then:
$$
\Lambda_{stat}({f_0}|_{\mathcal{B}(\Gamma_0)})=\bigcup_{a=1}^k \mathcal{P}_a\subset \Gamma_0.
$$
\end{proposition}

\textbf{Proof:}
The flow from $Out(\mathcal{P}_{a})$ to $In(\mathcal{P}_{a+1})$ is non-singular as in a flow-box. Since both $Out(\mathcal{P}_{a})$ and $In(\mathcal{P}_{a+1})$ are compact sets, the time of flight between them has a positive maximum. On the other hand, for each $a\in\ZZ_k$, the time of flight inside $V_a$ from $In(\mathcal{P}_{a})\backslash W^s_{loc}(\mathcal{P}_a)$ to $Out(\mathcal{P}_{a})$ tends to infinity as $t$ approach the stable manifold of $\mathcal{P}_a$,  $W^s_{loc}(\mathcal{P}_a)$, or equivalently as the trajectory accumulates on  $\Gamma_0$.
\Qed

 \begin{remark}\label{rkFlightTimes}
  It follows from Proposition \ref{density} that,  for each $a\in\ZZ_k$, the time intervals in which trajectories are travelling from $Out(\mathcal{P}_a)$ to $In(\mathcal{P}_{a+1})$ do not affect the accumulation points of the time averages of a solution that is accumulating on $\Gamma_0$. This result  will be useful in the proof of the Theorem \ref{Main1} because it shows that the duration of the journeys  between  nodes may be statistically neglected.
 \end{remark}

\subsection{Estimates of flight times}
In this section, we obtain relations between flight times of a trajectory in consecutive isolating blocks as well as other estimates that will be used in the sequel.
\begin{lemma}
\label{Lemma_times_2}
For all $j\in\NN$ and any initial condition $X\in  \mathcal{B}(\Gamma_0)$ we have:
\begin{equation}
\label{ratio1}
\frac{\tau_{j+1 }(X)}{\tau_{j}(X)}=\frac{c_{\langle j\rangle}}{e_{\langle j+1\rangle}}.
\end{equation}
In particular the ratio ${\tau_{j+1 }(X)}/{\tau_{j}(X)}$ does not depend on $X$. 
\end{lemma}

\textbf{Proof:}
Given  $j\in\NN$, let $X_j=\left(\theta_j,z_j\right)=\phi \left(T_j(X),X\right)\in In(\mathcal{P}_{\langle j\rangle})$.
 Using the expressions \eqref{Time of Flight}, \eqref{local map} and the expression for $\Psi_{\langle j\rangle}^0$ in \eqref{Psi_def}, we have:
$$
\tau_{j+1}(X)=\frac{1}{e_{\langle j+1\rangle}}\ln \left( \frac{\varepsilon}{ \varepsilon  \left(\frac{z_j}{\varepsilon}\right)^{\delta_{\langle j\rangle}}} \right) =
\frac{1}{e_{\langle j+1\rangle}} \delta_{\langle j\rangle} \left[\ln (\varepsilon)-\ln (z_j)\right]
$$
Thus
$$
 \frac{\tau_{j+1}(X)}{\tau_j(X)} = 
\frac{\frac{1}{e_{\langle j+1\rangle}} \delta_{\langle j \rangle} \left[\ln (\varepsilon)-\ln (z_j)\right]}{\frac{1}{e_{\langle j\rangle}}\left[\ln(\varepsilon)-\ln(z_j)\right]}
= \frac{e_{\langle j \rangle}}{e_{\langle j+1\rangle}} \delta_{\langle j\rangle} =\frac{c_{\langle j\rangle}}{e_{\langle j+1\rangle}}.
$$
\Qed

Recall  from \eqref{constants} that $\displaystyle\mu_{a+1}= \frac{c_a}{e_{a+1}}$,  $a \in \ZZ_k$.
With this notation we obtain:

\begin{corollary}
\label{Lemma3}
For $i,j\in \NN$ such that $j>i>1$,  and for any   $X\in  \mathcal{B}(\Gamma_0)$, we have:
\begin{enumerate}
\item 
\label{geom_sum}
$\frac{\tau_{j+k}(X)}{\tau_j (X)}=  \prod_{a=1}^k \mu_{a+1} =  \prod_{j=1}^k \delta_j=\delta>1$.
\item 
$\tau_{j+1}(X)= \tau_i (X)\prod^{j+1}_{l=i+1}\mu_{\langle l\rangle}$.
\end{enumerate}
\end{corollary}

We finish this section with a result comparing the two sequences of times $(T_i)_{i\in \NN}$ and $(\tau_i)_{i\in \NN}$.
The proof is very technical and is given in Appendix \ref{appendixB1}.

\begin{lemma} 
\label{Equalities}
 For $a\in\ZZ_k$, and for any  $X\in  \mathcal{B}(\Gamma_0)$, the following equalities hold:
\begin{enumerate}
\item $T_{a + nk}(X)=T_a (X)+ \frac{\delta^n-1}{\delta-1} \left(\mu_a +\mu_a \mu_{a+1} + \ldots + \prod_{l=0}^{k-1}\mu_{a+l}\right)\tau_{a-1}(X);$
\item $\tau_{a + nk}(X)=T_{a +1+ nk }(X)-T_{a + nk}(X)= \delta^n \mu_a \tau_{ a-1}(X)$.
\end{enumerate}
\end{lemma}

\subsection{The vertices of the $k$-polygon}\label{subsecvertices}
In this section, we show that in $\mathcal{B}(\Gamma_0)$ the time averages fail to converge, by finding several accumulation points for them.
For each $a\in\ZZ_k$, define the point 
\begin{equation}
\label{point1}
A_{a}= \frac{\overline{x}_{a} + \mu_{a+1} \overline{x}_{a+1} + \mu_{a+1}\mu_{a+2}  \overline{x}_{a+2} + \ldots+ \prod_{l=1}^{k-1}\mu_{a+l}  \overline{x}_{a+k-1}}{1+\mu_{a+1} +\mu_{a+1}\mu_{a+2}+\ldots+  \prod_{l=1}^{k-1}\mu_{a+l}}= \frac{num(A_a)}{den(A_a)}
\end{equation}
Note that $A_a$ and $num(A_a)$ lie in $\RR^3$ and $den(A_a)\in \RR$. Later we will see that these points are the vertices of a polygon of accumulation points. First we show that they are accumulation points for the time averages.

\begin{proposition}
\label{Prop6}
Let $a\in\ZZ_k$,  let $f_0$  be a vector field  in $\RR^3$  satisfying (P\ref{P1})--(P\ref{P3}) and let $\phi(t,X)$ a solution of $\dot{x}=f_0(x)$ with  $X\in\mathcal{B}(\Gamma_0)$. Then
$$
\lim_{n \rightarrow +\infty } \left[\frac{1}{T_{a+nk}} \int_0^{T_{a+nk}} \phi(t,X) dt \right]  = A_{a}
$$
\end{proposition}

In order to prove Proposition \ref{Prop6}, first we show that it is sufficient to consider the limit when $n\to\infty$ of the averages over one turn around $\Gamma_0$ and then we prove that these averages tend to $A_a$. The proof is divided in two technical lemmas, which may be found in Appendix \ref{appendixC}.

\subsection{The sides of the $k$-polygon}
In Section~\ref{subsecvertices} we have shown that for  $a \in \ZZ_k$, the time average over the sequences $T_{a+nk}$ of times accumulate, as $n\to\infty$ in the $A_a$. 
In this section we describe  accumulation points for intermediate sequences of times $t_n$.
For this, it will  be useful to know how $A_a$ and $A_{a+1}$ are related:

\begin{lemma}
\label{Colinear}
For all  $a \in \ZZ_k$, the following equalities hold:
\begin{equation}
\label{colinear2}
\mu_{a+1} den(A_{a+1}) =den (A_a)-(1-\delta) 
\qquad\mbox{and}\qquad
\mu_{a+1} num(A_{a+1})=num(A_a) - (1-\delta)\overline{x}_a .
\end{equation}
\end{lemma}

\begin{proposition}
\label{propColinear}
The point  $A_{a+1}$ lies in the segment connecting  $A_a$ to $\overline{x}_{a}$.
\end{proposition}

\textbf{Proof:}
We use Lemma~\ref{Colinear} to obtain
$$
num(A_{a+1}) =\frac{1}{ \mu_{a+1}}num(A_a)+\frac{\delta-1}{ \mu_{a+1}}\overline{x}_{a}
$$
and hence
$$
A_{a+1}=\frac{ num(A_{a+1}) }{den(A_{a+1}) }
=\left(\frac{den(A_{a}) }{ \mu_{a+1}den(A_{a+1}) }\right)\frac{num(A_a)}{den(A_{a}) }
+\left(\frac{\delta-1}{ \mu_{a+1}den(A_{a+1}) }\right)\overline{x}_{a}=\alpha A_a+\beta \overline{x}_{a}.
$$
Again from Lemma~\ref{Colinear} we have $den (A_a)=\mu_{a+1} den(A_{a+1})  -(\delta-1)$, and therefore
$$
\alpha=\frac{den (A_a)}{\mu_{a+1} den(A_{a+1}) }=1-\frac{\delta-1}{\mu_{a+1} den(A_{a+1}) }=1-\beta
$$
hence  $A_{a+1}$ lies in the line through $A_a$ and $\overline{x}_{a}$.
From the expression in  Lemma~\ref{Colinear} it follows that $\mu_{a+1} den(A_{a+1})-den (A_a)=  \delta-1<0$, 
hence $0<\alpha<1$ and thus  $A_{a+1}$ lies in the segment from $A_a$ to $\overline{x}_{a}$, 
proving the result.
\Qed

\begin{figure}
\begin{center}
\includegraphics[height=6cm]{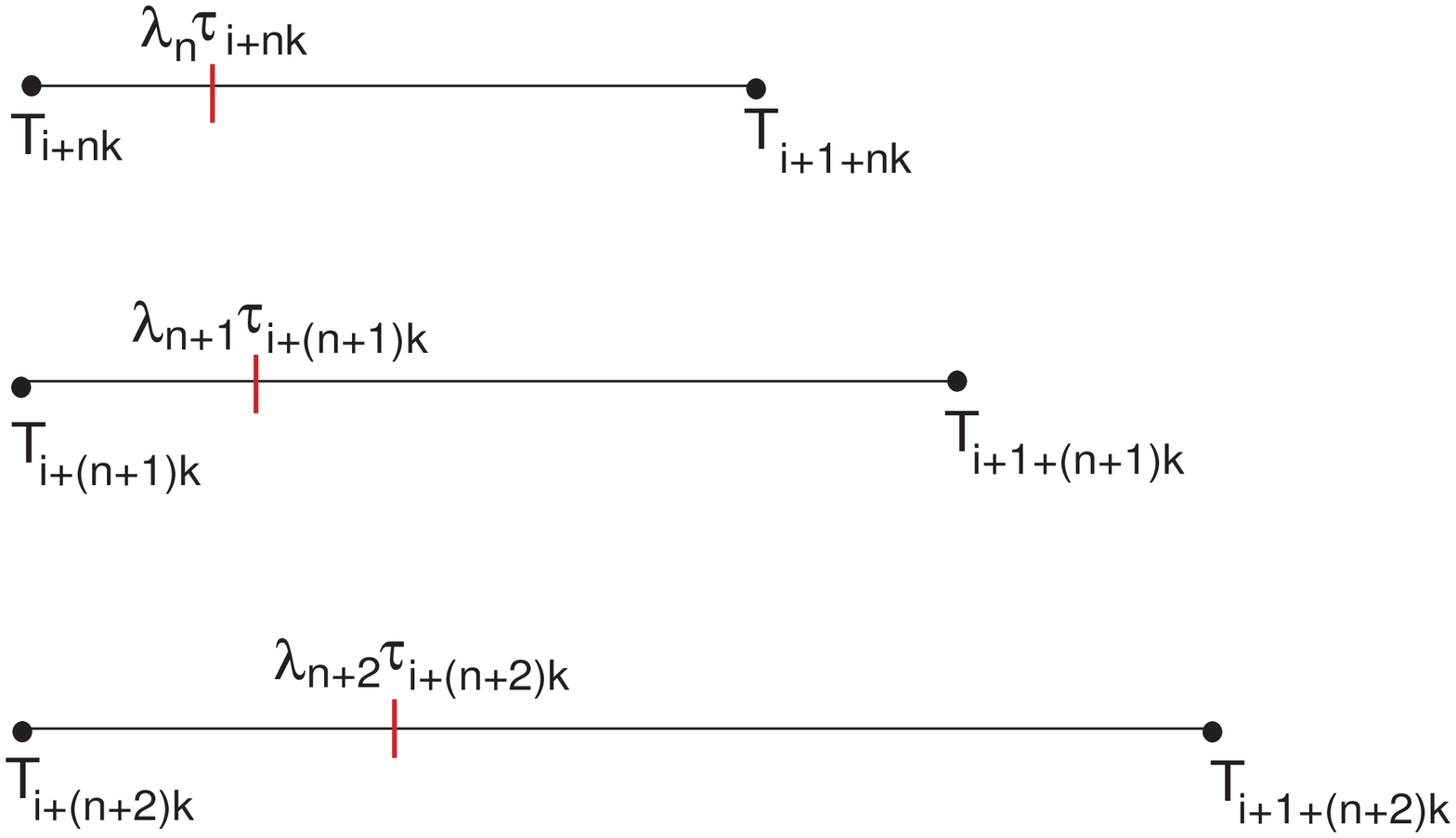}
\end{center}
\caption{\small Representation of the  sequence of times  $\lambda_n\tau_{a+nk}$, where $a\in\ZZ_k$,
is fixed and $n \in \NN$.}
\label{sequence2}
\end{figure}

We now come to the main result of this section:
\begin{theorem}
\label{Main1}
If $f_0$ is a vector field  in $\RR^3$ satisfying (P\ref{P1})--(P\ref{P3}), then
  for any $X\in  \mathcal{B}(\Gamma_0)$, the set of accumulation points of the time average $\frac{1}{T} \int_0^T \phi(t,X) dt $ is the boundary of the $k$-polygon defined by $A_1,\ldots, A_k\in \RR^3$  in (\ref{point1}). 
Moreover,
 when $\delta\rightarrow 1$ the polygon collapses into a point.
\end{theorem}

\textbf{Proof:}
First we show that all points in the boundary of the polygon are accumulation points.
Given $ L\in[0,1]$ and $a\in \ZZ_k$, consider the sequence $t_n=T_{a+nk}+ L \tau_{a+nk}$, we want the accumulation points of $\mathcal{L}_n=\frac{1}{t_n} \int_0^{t_n} \phi(t,X) dt$ as $n\to\infty$.
 For this we write
\begin{eqnarray*}
\mathcal{L}_n=
\frac{1}{t_n} \int_0^{t_n} \phi(t,X) dt  &=&
 \frac{1}{t_n} \int_0^{T_{a+nk}} \phi(t,X) dt + \frac{1}{t_n} \int_{T_{{a+nk}}}^{t_n} \phi(t,X) dt\\
&=&\alpha_n \left(\frac{1}{T_{{a+nk}}} \int_0^{T_{a+nk}} \phi(t,X) dt\right) 
+ \beta_n \left(\frac{1}{t_n-T_{{a+nk}}} \int_{T_{{a+nk}}}^{t_n} \phi(t,X) dt\right),
\end{eqnarray*}
where
$$
0< \alpha_n=\frac{T_{a+nk} }{t_n} \leq 1,
\qquad
0\leq \beta_n= \frac{t_n-T_{{a+nk}} }{t_n}\leq 1
\qquad\text{and} \qquad 
\alpha_n+\beta_n=1 .
$$
Since both $\alpha_n$ and $\beta_n$ are limited, each one of them contains a converging subsequence.
We analyse separately each of the terms in the expression for $\mathcal{L}_n$ above.

We have already seen in  Proposition~\ref{Prop6}   that,  if $X\in  \mathcal{B}(\Gamma_0)$,
then $\lim_{n\to\infty} \frac{1}{T_{a+nk}} \int_0^{T_{a+nk}} \phi(t,X) dt=A_a$.
In particular, if $ L=0$, then $\alpha_n=1$, $\beta_n=0$ and $\lim_{n\to\infty}\mathcal{L}_n =A_a$.

We claim that  if $ L\ne 0$, then
$\lim_{n\to\infty} \frac{1}{t_n-T_{a+nk}} \int_{T_{{a+nk}}}^{t_n} \phi(t,X) dt=\overline{x}_a$.
To see this, note that  $\phi(t,X)\in V_a$ for $t\in[T_{a+nk},t_n]$.
Moreover, since $\lim_{n\to\infty}\tau_{a+nk}=\infty$, then   for large $n$, we have that $t_n-T_{a+nk}= L\tau_{a+nk}$ is much larger than $\xi_a$, the period of $\mathcal{P}_a$.
Since $X\in  \mathcal{B}(\Gamma_0)$, then $ \phi(t,X) $, with  $t\in[T_{a+nk},t_n]$, tends to $\mathcal{P}_a$ when $n\to\infty$
and the average of $ \phi(t,X)$ tends to $\overline{x}_a$, the average of  $\mathcal{P}_a$.

At this point we have established that any accumulation point of $\mathcal{L}_n$ lies in the segment connecting $A_a$ to $\overline{x}_a$.
We have shown in Proposition~\ref{propColinear} that this segment  also contains $A_{a+1}$. 
By Proposition~\ref{Prop6} we have that $\lim_{n\to\infty}\mathcal{L}_n=A_{a+1}$ for $ L=1$.
On the other hand,  $\beta_n$ is an increasing function of $ L$, so, as $ L$ increases from 0 to 1, the accumulation points of $\mathcal{L}_n$ move from $A_a$ to $A_{a+1}$  in  the segment connecting them.

Conversely, any accumulation point lies  on the boundary of the polygon.
To see this,
let $A$ be an accumulation point of the time average.
This means that there is an increasing sequence of times $s_n$, tending to infinity, 
and such that 
$\lim_{n\to \infty}\mathcal{L}_n=A$, where $\mathcal{L}_n=\frac{1}{s_n}\int_0^{s_n}\phi(t,x)dt$.
Since $s_n$ tends to infinity, then it may be partitioned into subsequences of the form $s_{n_j}=T_{a+n_jk}+ L_{n_j} \tau_{a+n_jk}$ for each $a\in\ZZ_k$, and some $ L_{n_j}\in[0,1]$, as shown in Figure \ref{sequence2}.
The arguments above, applied to this subsequence, show that the accumulation points of $\mathcal{L}_{n_j}$ lie in the  segment connecting $A_a$ to $A_{a+1}$. 
Therefore, since $\mathcal{L}_n$ converges, there are two possibilities.
The first is that  all the $s_n$ (except possibly finitely many)  are of the form above for a fixed $a\in\ZZ_k$, and hence $A$ lies in the the  segment connecting $A_a$ to $A_{a+1}$.
The second possibility is that  all  the $s_n$ (except maybe a finite number)  are of one of the forms 
$$s_{n_j}=T_{a+n_jk}+ L_{n_j} \tau_{a+n_jk}\qquad  \text{or} \qquad s_{n_i}=T_{a+1+n_ik}+ L_{n_i} \tau_{a+1+n_ik},$$
 and that $A=A_{a+1}$.
In both cases, the accumulation point of the time average will lie on the boundary of the polygon. Finally, when $\delta\to 1$, the expressions \eqref{colinear2} in Lemma~\ref{Colinear} become
$\mu_{a+1} den(A_{a+1}) =den (A_a)$ and $\mu_{a+1} num(A_{a+1})=num(A_a)$, hence
$$
A_{a}=\frac{num(A_a)}{den (A_a)}=\frac{\mu_{a+1} num(A_{a+1})}{\mu_{a+1} den(A_{a+1}) }=A_{a+1}
$$
and the polygon collapses to a point at the same time as $\Gamma_0$ stops being attracting.
\Qed

\begin{figure}
\begin{center}
\includegraphics[height=6cm]{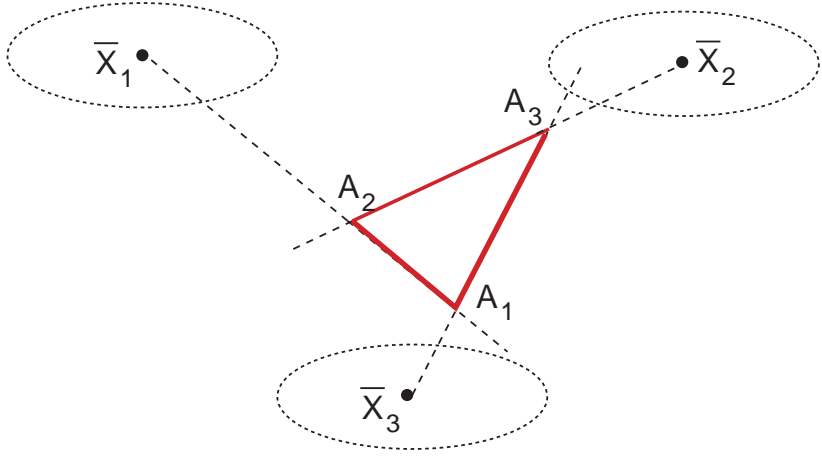}
\end{center}
\caption{\small The polygon in Theorem \ref{Main1} with $k=3$: the accumulation points of the time average $\frac{1}{T} \int_0^T \phi(t,x) dt $  lie on the boundary of the triangle defined by $A_1,A_2$ and $A_3.$ }
\label{scheme1}
\end{figure}

Taking the observable as the projection on any component, the first main result of this paper may be stated as:

\begin{corollary}
\label{CorollaryHistoric}
If $f_0$ is a vector field  in $\RR^3$ satisfying (P\ref{P1})--(P\ref{P3}), then 
all points in the basin of attraction of $\Gamma_0$ have historic behaviour. In particular the set of initial conditions with historic behaviour has positive Lebesgue measure.
\end{corollary}

The  points of  $\Gamma_0$ do not have historic behaviour. 
Indeed, if  $X\in\Gamma_0$ then either $X\in \mathcal{P}_a$ or $\phi(t,X)$ accumulates on $\mathcal{P}_a$ for some $a\in\{1,\ldots,k\}$.
In both cases, $\lim_{T\to\infty}\frac{1}{T}\int_0^T\phi(t,X)dt=\overline{x}_a$.
\medbreak
The previous proofs have been done for a piecewise continuous trajectory;  when $t=T_a$, the trajectory jumps from $V_{a-1}$ to $V_a$, whereas the real solutions have a continuous motion from $V_{a-1}$ to $V_a$ along the corresponding heteroclinic connection, during a bounded interval of time. 
 As shown in Proposition \ref{density}, the statistical limit set of $\Gamma_0$ is $\bigcup_{a=1}^k \mathcal{P}_a$  meaning that trajectories spend Lebesgue almost all time near the periodic solutions, and not along the connections. 
 Therefore, the intervals in which the transition occurs do not affect the accumulation points of the time averages of the trajectories and the result that was shown for a piecewise continuous trajectory holds.

\section{Persistence of  historic behaviour}
\label{Tangencies}
From now on, we discuss he differential equation $\dot{x}=f_\lambda(x)$ satisfying (P\ref{P1})--(P\ref{P5}), with 
$\lambda\ne 0$.
In this case it was shown  in Rodrigues \emph{et al} \cite{Rodrigues} that  the simple dynamics near $\Gamma_0$ jumps to chaotic behaviour near $\Gamma_\lambda$.

\subsection{Invariant manifolds for $\lambda>0$}
\label{secTransition}
\begin{figure}[htt]
\begin{center}
\includegraphics[height=4.5cm]{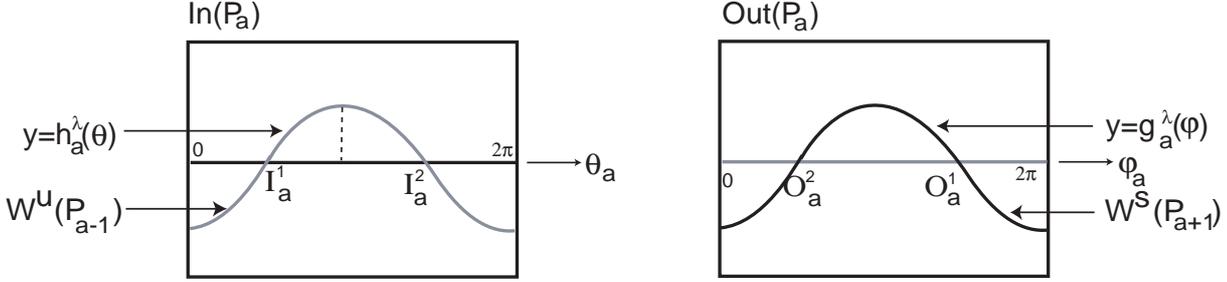}
\end{center}
\caption{\small For $\lambda$ close to zero, both $W^s_{loc}(\mathcal{P}_{a+1})\cap Out(\mathcal{P}_{a})$
and $W^u_{loc}(\mathcal{P}_{a})\cap In(\mathcal{P}_{a+1})$ are closed curves, given in local coordinates as the  graphs of  periodic functions; this is the expected unfolding from the coincidence of the invariant manifolds at $\lambda=0$.}
\label{Transitions}
\end{figure}
\medbreak
We describe the geometry of the two-dimensional local invariant manifolds of $\mathcal{P}_{a}$ and $\mathcal{P}_{a+1}$ for  $\lambda\neq 0$, under the assumptions  (P\ref{P1})--(P\ref{P5}).
For this, let $f_\lambda$ be an unfolding of $f_0$ satisfying (P\ref{P1})--(P\ref{P5}). 
For $\lambda\ne 0$, we introduce the notation:
\begin{itemize}
\item 
$(O_a^1,0)$ and $(O_a^2,0)$ with $0<O_a^1<O_a^2<2\pi$ are the coordinates of the two points
 where the connections  $[\mathcal{P}_{a} \rightarrow \mathcal{P}_{a+1}]$ of Properties (P\ref{P4})--(P\ref{P5})   meet $Out(\mathcal{P}_{a})$;
\item
$(I_a^1,0)$ and $(I_a^2,0)$  with $0<I_a^1<I_a^2<2\pi$  are the coordinates of the two  points  
where $[\mathcal{P}_{a-1} \rightarrow \mathcal{P}_a]$ meets $In(\mathcal{P}_{a})$;
\item
$(O_a^i,0)$ and $(I_{a+1}^i,0)$ are on the same trajectory for each $i \in \{1,2\}$ and $a\in\ZZ_k$.
\end{itemize}

By (P\ref{P5}), for small $\lambda>0$, the curves $W^s_{loc}(\mathcal{P}_{a+1})\cap Out(\mathcal{P}_a)$ and $W^u_{loc}(\mathcal{P}_{a})\cap In(\mathcal{P}_{a+1})$
can be seen as graphs of smooth periodic functions, for which we  make the following conventions (see Figure \ref{Transitions}):
\begin{itemize}
\item
$W^s_{loc}(\mathcal{P}_{a+1})\cap Out(\mathcal{P}_{a})$ 
is the graph of   $y=g_a^\lambda (\varphi)$, with $g_a^\lambda(O_a^i)=1$, for $i \in \{1,2\}$ and $a\in\ZZ_k$.
\item
$W^u_{loc}(\mathcal{P}_{a-1})\cap In(\mathcal{P}_{a})$ is the graph
  of  $y=h_{a}^\lambda (\theta)$, with $h_{a}^\lambda(I_{a}^i)=0$, for $i \in \{1,2\}$ and $a\in\ZZ_k$.
\item
omitting the superscript $\lambda$, we have: $h_a^\prime(I_a^1)>0$, $h_a^\prime(I_a^2)<0$, $g_a^\prime(O_a^2)>0$ and $g_a^\prime(O_a^1)<0$, for $i \in \{1,2\}$.
\end{itemize}

\begin{figure}[hbb]
\begin{center}
\includegraphics[height=4cm]{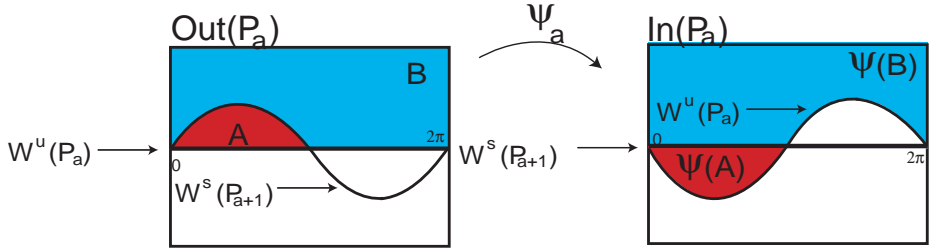}
\end{center}
\caption{\small The component $A$ of  $Out(\mathcal{P}_{a})$ between $W_{loc}^s(\mathcal{P}_{a+1})$ and  $W_{loc}^u(\mathcal{P}_{a})$ 
 is mapped by $\Psi_a$ into the lower ($z_{a+1}<0$) part of $In(\mathcal{P}_{a+1})$, its complement $B$  in the $r_a>1$ component of $Out(\mathcal{P}_{a})$  is mapped by $\Psi_a$ into the upper ($z_{a+1}>0$)  part of $In(\mathcal{P}_{a+1})$.} 
\label{figRegionsInOut}
\end{figure}

The two points $(O_a^1,0)$ and $(O_a^2,0)$ divide the closed curve $W^s_{loc}(\mathcal{P}_{a+1})\cap Out(\mathcal{P}_{a})$
  in two components, corresponding to different signs of $r_a-1$.
 With the conventions  above, we get $g_a^\lambda(\varphi)>1$ for  $\varphi \in\left(O_a^2,O_a^1\right)$.
 More specifically,  the region in $Out(\mathcal{P}_{a})$ between $W_{loc}^s(\mathcal{P}_{a+1})$ and  $W_{loc}^u(\mathcal{P}_{a})$ given by 
 $$
 A=\{(\varphi_a,r_a)\in Out(\mathcal{P}_{a}): 1<r_a<g_a^\lambda(\varphi_a) \}
 $$
is mapped by $\Psi_a$ into the lower ($z_a<0$) part of $In(\mathcal{P}_{a+1})$.
Similarly, the region 
$$
B=\{(\varphi_a,r_a)\in Out(\mathcal{P}_{a}): r_a>1\}\backslash A=\{(\varphi_a,r_a)\in Out(\mathcal{P}_{a}): 1<r_a\ \mbox{ and }\ g_a^\lambda(\varphi_a)<r_a \}
$$ 
(see Figure~\ref{figRegionsInOut})
is mapped into the $z_a>0$ component of $In(\mathcal{P}_{a+1})$.

The maximum value of $g_a^\lambda (\varphi)$ is  attained at some point 
$$
(\varphi_a,r_a)= (\varphi_a^O (\lambda),M_a^O(\lambda)) \qquad \text{with} \qquad O_a^2<\varphi_a^O(\lambda)<O_a^1.
$$
We denote by $M^I_a(\lambda)$ the maximum value of $h_a^\lambda$.

\subsection{Geometrical preliminaries}
We will need to  introduce some  definitions. 

\begin{figure}
\begin{center}
\includegraphics[height=6cm]{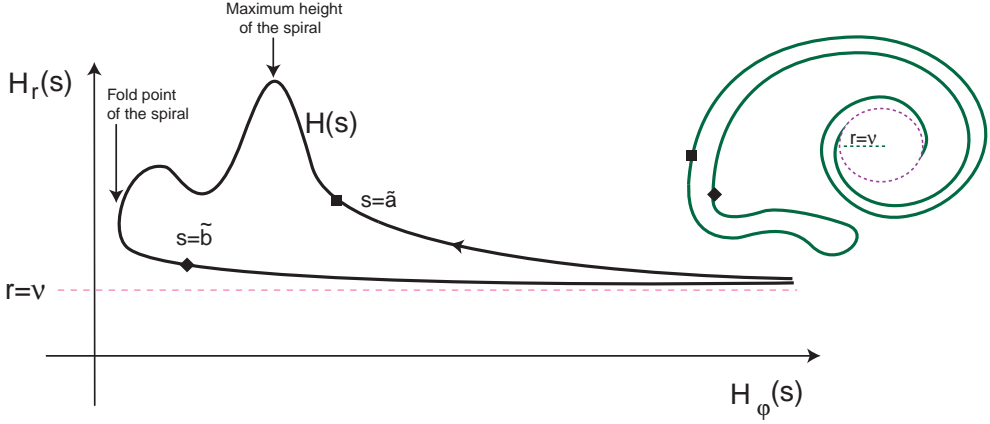}
\end{center}
\caption{\small A  spiral is defined on a covering of the annulus $Out(\mathcal{P}_a)$ by a smooth curve that turns around the annulus infinitely many times as its radius tends to $\nu\in[0,1]$. It contains a fold point and a point of maximum radius. }
\label{spiral1}
\end{figure}

\begin{definition}
\label{spiral_def}
A \emph{spiral} on the annulus $\mathcal{A}$ \emph{accumulating on the circle} $r=\nu$ is a curve on $\mathcal{A}$, without self-intersections, that is the image, by the parametrisation $(\varphi,r )$ of the annulus, of a continuous map
$H:(b,c)\rightarrow \RR\times[0,1]$, 
$$
H(s)=\left(\varphi(s),r(s)\right),
$$
 such that:
\begin{enumerate}
\renewcommand{\theenumi}{\roman{enumi}}
\renewcommand{\labelenumi}{{\theenumi})}
\item \label{monotonicity}
there are $\tilde{b}\le \tilde{c}\in (b,c)$ for which both $\varphi(s)$ and $r(s)$ are monotonic in each of the intervals $(b,\tilde{b})$ and $(\tilde{c},c)$;
\item \label{turns}
either $\lim_{s\to b^+}\varphi(s)=\lim_{s\to c^-}\varphi(s)=+\infty$ or $\lim_{s\to b^+}\varphi(s)=\lim_{s\to c^-}\varphi(s)=-\infty,$
\item \label{accumulates}
$\lim_{s\to b^+}r(s)=\lim_{s\to c^-}r(s)=\nu$.
\end{enumerate}
\end{definition}

It follows from the assumptions on the function $\varphi(s)$ that it has either a global minimum or a global maximum, and that $r(s)$ always has a global maximum.
The point where the map $\varphi(s)$  has a global minimum or a global maximum will be called a \emph{fold point} of the spiral.
The global maximum value of $r(s)$ will be called the \emph{maximum radius} of the  spiral.

\subsection{Geometry of the transition maps $\Phi_a$}
\label{secImageInvariantManifs}

\begin{proposition}
\label{Structures} Under the conventions of Section~\ref{secTransition}, 
for each $a \in \ZZ_k$, the local map $\Phi_a$ transforms the part of the graph of $h_a$ with $I_a^1<\theta<I_a^2$ into a spiral on $Out(\mathcal{P}_{a})$ accumulating on the circle  $Out(\mathcal{P}_{a}) \cap W^{u}_{loc}(\mathcal{P}_{a})$.
This spiral  has maximum radius $1+ \varepsilon^{1-\delta_a}(M_a^I)^{\delta_a}$;
it has a fold point that, as $\lambda$ tends to zero, turns around $Out(\mathcal{P}_{a})$ infinitely many times.
\end{proposition}

 \textbf{Proof:}
 The curve 
  $\Phi_a\left(W^u_{loc}(\mathcal{P}_{a-1})\cap In(\mathcal {P}_a  \right)$ is given by 
$H_a(\theta)= \Phi_a(\theta, h_a(\theta)) =(\varphi_a(\theta), r_a(\theta))$ where:
\begin{equation}
\label{spiral_expression}
H_a(\theta)= \Phi_a(\theta, h_a(\theta)) = \left(\theta-\frac{1}{e_a}\ln\left(\frac{h_a(\theta)}{\varepsilon}\right),
1+ \varepsilon \left(\frac{h_a(\theta)}{\varepsilon}\right)^{\delta_a}\right)=
(\varphi_a(\theta), r_a(\theta)) .
\end{equation}
From this expression if follows immediately that 
$$
\lim_{\theta \rightarrow I_a^1} \varphi_a(\theta)=\lim_{x \rightarrow I_a^2} \varphi_a(\theta)= +\infty 
\quad\mbox{and}\quad
\lim_{\theta \rightarrow I_a^1} r_a(\theta)=\lim_{\theta \rightarrow I_a^2} r_a(\theta)= 1
$$
hence, conditions {\sl \ref{turns})} and {\sl \ref{accumulates})} of the definition of  spiral hold.
Condition {\sl \ref{monotonicity})} holds trivially near $I_a^2$ since $h_a'(I_a^2) < 0$, hence there is  
$\tilde{I_a^2}< I_a^2$ such that $\varphi_a'(\theta)>1$ for  all $\theta \in \left( \tilde{I_a^2}, I_a^2\right)$.
On the other hand,  since $h_a'(I_a^1) > 0$ and $\lim_{\theta \rightarrow I_a^1}h_a(\theta)=0$, there is  $\tilde{I_a^1}<\theta_a^M$, where $\varphi_a'(\theta)<0$ for all $\theta \in \left(I_a^1, \tilde{I_a^1}\right)$.

 The statement about the maximum radius follows immediately from \eqref{spiral_expression}  and the conventions of Section~\ref{secTransition}.
 
Let $H_a(\theta_a^\star(\lambda))$ be a fold point of the spiral. 
Its first coordinate is given by 
$\varphi_{a}^\star=\theta_a^\star-\frac{1}{e_a}\ln\left(\frac{h_a(\theta_a^\star(\lambda))}{\varepsilon}\right)$ and $h_a(\theta_a)\le M^I_a(\lambda)$.
Since $f_\lambda$ unfolds $f_0$, then $\lim_{\lambda\to 0}M^I_a(\lambda)=0$ and therefore
$\lim_{\lambda\to 0}\varphi_{a}^\star=+\infty$.
 Hence, the fold point  turns around the cylinder $Out(\mathcal{P}_{a})$  infinitely many times, as $\lambda$ tends to zero. 
\Qed

\subsection{A set of one-parameter families of  vector fields}

For any  unfolding $f_\lambda$ of $f_0$, as we have seen in Sections~\ref{secTransition} and \ref{secImageInvariantManifs},  the maximum radius $M_a^O(\lambda)$ of $W^s_{loc}(\mathcal{P}_{a+1})\cap Out(\mathcal{P}_{a})$, and the maximum height $M_a^I(\lambda)$ of $W^u_{loc}(\mathcal{P}_{a-1})\cap In(\mathcal{P}_{a})$, satisfy:
$$
\lim_{\lambda\to 0} M_a^I(\lambda)=0\qquad
\lim_{\lambda\to 0} \left(1+ \varepsilon^{1-\delta_a}(M_a^I(\lambda))^{\delta_a}\right)=\lim_{\lambda\to 0} M_a^O(\lambda)=1.
$$

We make the additional assumption that  $1+ \varepsilon^{1-\delta_a}(M_a^I(\lambda))^{\delta_a}$ tends to zero faster than $M_a^O(\lambda)$ for at least one $a \in\ZZ_k$. 
 This condition defines the open set ${\mathcal C}$ of generic unfoldings $f_\lambda$ that we need for the statement of Theorem \ref{teorema tangency}.  
 More precisely,
 \begin{equation}
 \label{eqDefineC}
 {\mathcal C}=\left\{ 
 f_\lambda \mbox{ satisfying (P\ref{P1}) -- (P\ref{P5})}: \exists a\in\ZZ_k\ 
 \exists \lambda_0>0 :\ \ 
 0<\lambda<\lambda_0\ \Rightarrow
1+ \varepsilon^{1-\delta_a}(M_a^I(\lambda))^{\delta_a}<M_a^O(\lambda)
 \right\}.
 \end{equation}
 The   set ${\mathcal C}$ is open in the Whitney  $C^2$ topology. 

\subsection{Heteroclinic tangencies}

\begin{theorem}
\label{teorema tangency}
For any family  $f_\lambda$ of vector fields in the set ${\mathcal C}$ defined in \eqref{eqDefineC}  there is $a\in\ZZ_k$ such that:
\begin{enumerate}
\item\label{tangentManifs}
 there is a sequence $\lambda_i>0$ of real numbers with $\lim_{i\to\infty}\lambda_i=0$ such that
for $\lambda=\lambda_i$ the manifolds $W^u(\mathcal{P}_{a-1})$ and $W^s(\mathcal{P}_{a+1})$ are tangent;
for $\lambda>\lambda_i$,
there are two heteroclinic connections in $W^u(\mathcal{P}_{a-1})\cap W^s(\mathcal{P}_{a+1})$ that  collapse into the tangency at $\lambda=\lambda_i$ and then disappear for $\lambda<\lambda_i$; 
\item\label{accumulatingTangs}
 arbitrarily close to the connection $[\mathcal{P}_{a-1}\to \mathcal{P}_{a}]$
there are  hyperbolic periodic   solutions at points $x_i$  and infinitely many values $\lambda_{n,i}$ 
for which  the periodic  solution has a homoclinic tangency of its  invariant manifolds.
\end{enumerate}

\end{theorem}

\begin{figure}
\begin{center}
\includegraphics[height=4cm]{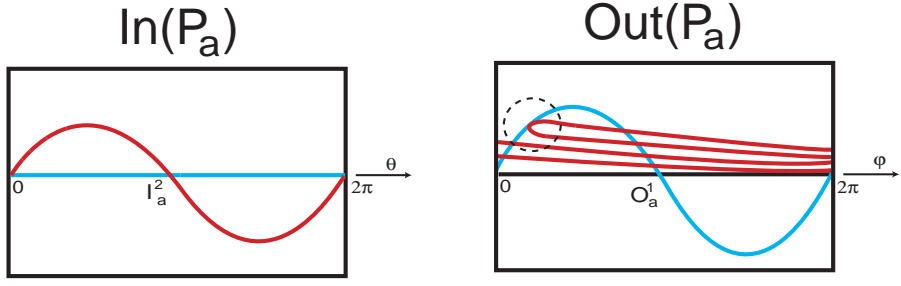}
\end{center}
\caption{\small When $\lambda$ decreases, the fold point of the spiral $\Phi_a\left(W^u_{loc}(\mathcal{P}_{a-1})\cap In(\mathcal {P}_a)  \right)$
 moves to the right and for $\lambda=\lambda_i$, it is tangent to $W^s(\mathcal{P}_{a+1})$ creating a heteroclinic tangency.  }
\label{tangencies1}
\end{figure}

Note that for $k=2$, the tangency  of assertion {\sl \ref{tangentManifs}.} is a homoclinic connection.

 \textbf{Proof:}
 Let $\theta_a=\theta_a^\star(\lambda)$ correspond to a  fold point of the  
 spiral   $\Phi_a\left(W^u_{loc}(\mathcal{P}_{a-1})\cap In(\mathcal {P}_a)  \right)$ given by \eqref{spiral_expression}.
Since $f_\lambda\in{\mathcal C}$ and using Proposition~\ref{Structures} and \eqref{eqDefineC},  for $\lambda< \lambda_0$  all points in the 
spiral have second coordinate less than $M_a^O$, this is true, in particular, for  the fold point $H_a(\theta_a^\star(\lambda))$.
Also by Proposition~\ref{Structures} the  fold point turns around $Out(\mathcal{P}_{a})$ infinitely many times as $\lambda$ goes to zero.
This means that there is a positive value $\lambda_A<\lambda_0$ such that 
$H_a(\theta_a^\star(\lambda_R))$ lies in the region $A$  that will be mapped to $z_a<0$ (see   Section~\ref{secTransition}) and there is  a positive value $\lambda_B<\lambda_A$ such that $H_a(\theta_a({\lambda_L}))$ lies in the region $B$ that goes to $z_a>0$, as in Figure~\ref{tangencies1}.
Therefore,  the curve $H_a(\theta_j(\lambda))$ is tangent to the graph of $g_a^\lambda$ at some point $H_a(\theta_a^\star(\lambda_1))$ with $\lambda_1\in\left( \lambda_B,\lambda_A\right)$.

As $\lambda$ decreases from $\lambda_B$, the fold point enters and leaves the region $A$, creating a sequence of tangencies  to the graph of $g_a^\lambda$.  
At each tangency, two points where $H_a(\theta_a^\star(\lambda))$ intersects the graph of   $g_a^\lambda$ come together, corresponding to the pair of transverse heteroclinic connections that collapse at the tangency.
This completes the proof of {\sl \ref{tangentManifs}.}

\begin{figure}
\begin{center}
\includegraphics[width=11cm]{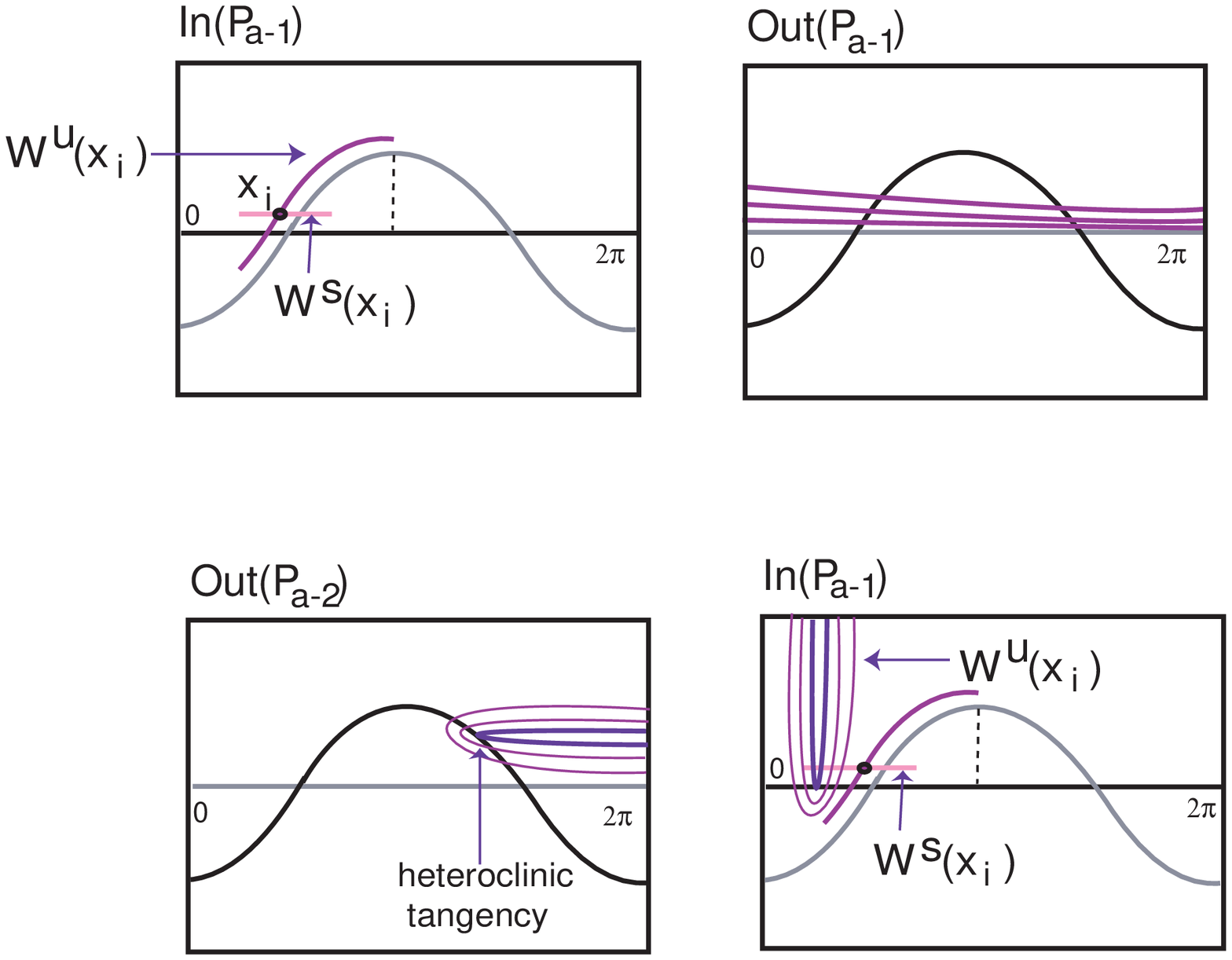}
\end{center}
\caption{\small The unstable manifold of a fixed point $x_i$ of the first return map to $In(\mathcal{P}_{a-1})$  accumulates on $W^u(\mathcal{P}_{a-1})$ and defines a family of curves in $Out(\mathcal{P}_{a})$ with  a fold point.
When $\lambda$ decreases, the fold point  moves to the right and for $\lambda=\lambda_{n,i}$, it is tangent to $W^s(x_i)$ creating a homoclinic tangency.  }
\label{homoclinicTangency}
\end{figure}

For assertion   {\sl \ref{accumulatingTangs}.}, note that by the results of \cite{Rodrigues} there is a suspended horseshoe near the connection $[\mathcal{P}_{a-1}\to \mathcal{P}_{a}]$.
Hence, there are hyperbolic fixed points of the first return map to $In(\mathcal{P}_{a-1})$ arbitrarily close to the connection; let  $p_i$ be one of them.
Denote by  $\eta_a$ the map  $\Psi_{a } \circ \Phi_a$. 
The image by  $\Phi_{a-1}$ of an interval contained in $W^u(p_i)$ accumulates on $W^u(\mathcal{P}_{a-1})$. 
In
particular, it is mapped by $\eta_a\circ\Phi_{a-1}$
into  infinitely many spirals in $Out(\mathcal{P}_{a})$, each one having  a fold point --- see Figure~\ref{homoclinicTangency}.
Since the fold points turn around $Out(\mathcal{P}_{a+1})$ infinitely many times as $\lambda$ varies,
this curve is tangent to $W^s(p_i)$ at a sequence $\lambda_{n,i}$ of values of $\lambda$.
\Qed

The hypothesis (P\ref{P3}) in the definition of $\mathcal C$ for Theorem~\ref{teorema tangency}, that the family $f_\lambda$ unfolds the  degeneracy $f_0$, may be replaced by the assumption that the flow of $f_\lambda$ turns in opposite directions around two successive nodes $\mathcal{P}_{a}$ and $\mathcal{P}_{a+1}$, as in \cite{LR3}, {\sl ie}, by the assumption that two successive nodes have different chirality.
This is because, in the proof of Theorem~\ref{teorema tangency}, the heteroclinic tangency is obtained from the presence of a fold pont in the curve $\Phi_a\left(W^u_{loc}(\mathcal{P}_{a-1})\cap In(\mathcal {P}_a)  \right)$ and from the control of the angular coordinate $\varphi$ of the fold. This is the content of Proposition~\ref{Structures}, where we use the fact that $W^u_{loc}(\mathcal{P}_{a-1})\cap In(\mathcal {P}_a) $ is the graph of a function with a maximum, a consequence of 
 (P\ref{P3}) and  (P\ref{P4}). 
 If we assume instead that successive nodes have different chirality as in \cite{LR3}, then  the image by $\Phi_a$ of the curve $W^u_{loc}(\mathcal{P}_{a-1})\cap In(\mathcal {P}_a) $ will have infinitely many fold points  whose coordinates $\varphi$ will form a dense subset of $[0,2\pi]$, and hence, as in \cite{LR3}, an arbitrarily small change in the parameter $\lambda$ will create a heteroclinic tangency.

\subsection{Historic behaviour}
 The next result is the core of this section.
It locates trajectories with historic behaviour  $C^2$-close to the unfolding of a degenerate equation, as a consequence of the tangencies found in Theorem~\ref{teorema tangency}.

\begin{theorem}
\label{teoremaHistoric}
For any family  $f_\lambda$ of vector fields in the open set ${\mathcal C}$ defined in \eqref{eqDefineC} 
there are sequences $0<\xi_i<\zeta_i<\xi_{i+1}$, with $\lim_{i\rightarrow +\infty}\zeta_i=0$,
such that for each $\lambda$ in $(\xi_i,\zeta_i)$, there are vector fields arbitrarily close to $f_\lambda$ in the $C^2$-topology for which there is an open set of initial conditions with historic behaviour.
\end{theorem}

 In the proof we will use the following concept:

\begin{definition}
Let $M$ be a smooth surface and  let $\mathrm{Diff}^r(M)$  be the set of its local diffeomorphisms of class $C^r$, $r\ge 2$.
An open subset $\mathcal{N}\subset \mathrm{Diff}^r(M)$ is a \emph{Newhouse domain} if any element of $\mathcal{N}$ is $C^r$-approximated  by a diffeomorphism $g$ with a homoclinic tangency associated with a dissipative saddle fixed point $p_g$, and moreover $g$ has a $C^r$-persistent tangency associated with some basic sets $\Lambda_g$ containing $p_g$ in the sense that there is a $C^r$-neighbourhood of $g$ any element of which has a homoclinic tangency for the continuation of $\Lambda_g$.
\end{definition}

Newhouse has shown in \cite{Newhouse79} that any $C^2$ diffeomorphism containing a homoclinic tangency to a dissipative saddle point lies in the closure of a Newhouse domain in the $C^2$ topology.

We will need the definition of historic behaviour for diffeomorphisms:
\begin{definition}
Let $F$ be a $C^2$ diffeomorphism on a smooth surface $M$.
We say that the forward orbit 
$
\{ x,F(x),F^2(x),\ldots, F^j(x),\ldots\}
$
has \emph{historic behaviour}
if the average
\begin{equation}\label{discreteHistoric}
\frac{1}{n+1}\sum_{j=0}^n \delta_{F^j(x)}
\end{equation}
does not converge as $n\to+\infty$ in the weak topology, where $\delta_Z$ is the Dirac measure on $M$ supported at $Z\in M$.
\end{definition}

 \textbf{Proof of Theorem~\ref{teoremaHistoric}:}
 For $\lambda=0$ and $a\in \ZZ_k$, the derivative of  the first return map  to $In(\mathcal{P}_{a})$ has determinant of the form $Cz_a^{\delta-1}$  for some constant $C>0$.
Thus, for sufficiently small $\lambda>0$, and at points near  $W^s(\mathcal{P}_{a})$, the first return map to $In(\mathcal{P}_{a})$  is also contracting, since the determinant  of its derivative has absolute value less than 1.
 Moreover, the family $f_\lambda$ unfolds  each one of the  homoclinic tangencies of Theorem~\ref{teorema tangency} generically.
Hence the arguments of Newhouse, Palis \& Takens and Yorke \& Alligood \cite{Newhouse79,PT, YA} revived in \cite{LR2015}  may be applied here to show that near each one of the  homoclinic tangencies  there is a sequence of intervals $(\xi_i,\zeta_i)$ in the set of parameters $\lambda$ corresponding to a Newhouse domain.

By Theorem A of Kiriki \& Soma \cite{KS}, each Newhouse domain for the first return map is contained in the closure of the set of diffeomorphisms having an open set of points with historic behaviour. 
 By Theorem~\ref{teorema tangency} the family $f_\lambda$ unfolds the heteroclinic tangencies generically.
 Hence, from the results of  \cite{KS},  it follows that for each $\lambda\in(\xi_i,\zeta_i)$, the first return map $F_a^\lambda$ may be approximated in the $C^2$ topology by maps $\widehat{F}$ defined in $In(\mathcal{P}_{a})$ for which we may find an open connected subset $\mathcal{U}\subset In(\mathcal{P}_{a})$ and two sequences of integers,  $(a_j)_{j\in \NN}$ and $(b_j)_{j\in \NN}$,  such that, for each $x\in\mathcal{U}$, the limits \eqref{discreteHistoric} for $\widehat{F}$ are different for $n$ in the two sequences.
In particular, there exists a set $A\subset In(\mathcal{P}_{a})$, such that 
$$
\forall x \in \mathcal{U}, \qquad 
\lim_{k \rightarrow +\infty} \frac{1}{a_k+1}\sum_{j=0}^{a_k} \delta_{\widehat{F}^j(x)}(A)
\neq \lim_{k \rightarrow +\infty} \frac{1}{b_k+1}\sum_{j=0}^{b_k} \delta_{\widehat{F}^j(x)}(A)
$$
or, equivalently,
\begin{equation}
 \label{Hist_Beha1}
\forall x \in \mathcal{U}, \qquad 
L=\lim_{k \rightarrow +\infty} \frac{1}{a_k+1}\sum_{j=0}^{a_k} \chi_A(\widehat{F}^{j}(x)) 
\neq \lim_{k \rightarrow +\infty} \frac{1}{b_k+1}\sum_{j=0}^{b_k} \chi_A(\widehat{F}^{j}(x)) = \widehat{L}
 \end{equation}
 where $\chi_A$ denotes the characteristic function on $A$.
 According to the proof of \cite{KS}, the two fixed points of the horseshoe  that arises near the tangency will be visited by orbits of points in the set 
 $\mathcal{U}$.

Since the maps $\widehat{F}$ are close to the first return map in the $C^2$ topology, they
may be seen as the first return maps to a vector field $g$ that is $C^2$-close to $f_\lambda$ --- see, for instance Remark 2 in Pugh and Robinson \cite[Section 7A]{PR}.
It remains to show that solutions to  $\dot{x}=g(x)$ have historic behaviour in the sense of Definition~\ref{historicDef}.

Let $\tau(x)$ be the time of first return of $x\in In (\mathcal{P}_a)$, \emph{ie} $\tau(x)>0$ and $\phi(\tau(x),x)\in In (\mathcal{P}_a)$ where $\phi$ is the flow associated to $\dot{x}=g(x)$.
Since $\mathcal{U}$ is connected, taking its closure $\overline{\mathcal{U}}$ compact and sufficiently small, 
$\tau(x)$ is approximately constant on $\overline{\mathcal{U}}$. 
Rescaling the time $t$ we may suppose $\tau(x)\equiv 1$.

Given $b>0$, let $V_b=\{\phi(t,x):\ -b<t<b, \ x\in In(\mathcal{P}_a)\}$.
For $0<c<1$ and $\varepsilon>0$ sufficiently small,  let $\psi: \RR^3 \rightarrow [0,1]$ be of class $C^k$, $k\ge 2$, such that
$\psi=$1 on $\overline{V_c}$ and $\psi=0$ outside $V_{c+\varepsilon}$.
Let $\mathcal{S}(A)$ be the saturation of $A$ by the flow $\phi$, given by
$\mathcal{S}(A)=\{\phi(t,x):\ t\in\RR, \ x\in A\}$.
Define the observable $H$ by $H(x)=\psi(x) \chi_{\mathcal{S}(A)}(x)$.
For $x\in\mathcal{U}$ we have:
\begin{eqnarray*}
\int_0^{a_{k}+1}H(\phi(t,x)) dt &=& \int_{0}^{c} H(\phi(t,x)) dt  +  \sum_{j=1}^{a_k}  \int_{j-c}^{j+c} H(\phi(t,x)) dt + \int_{a_{k}+1-c}^{a_{k}+1} H(\phi(t,x)) dt +o(\varepsilon) \\
&=&c  \chi_A(\widehat{F}^{0}(x)) +2c\sum_{j=1}^{a_k}   \chi_A(\widehat{F}^{j}(x)) +c  \chi_A(\widehat{F}^{a_{k}+1}(x)) +o(\varepsilon) .
\end{eqnarray*}
Hence
$$
\frac{1}{a_{k}+1}\int_0^{a_{k}+1}H(\phi(t,x)) dt =
\frac{2c}{a_{k}+1}\sum_{j=0}^{a_k}   \chi_A(\widehat{F}^{j}(x)) 
- \frac{c}{a_{k}+1} \chi_A(x)+\frac{c}{a_{k}+1} \chi_A(\widehat{F}^{a_{k}+1}(x))
$$
where
$$
\lim_{k \to +\infty}  \frac{c}{a_{k}+1}\chi_A(x)=0
\quad\mbox{and}\quad
\lim_{k\to +\infty} \frac{c}{a_{k}+1} \chi_A(\widehat{F}^{a_{k}+1}(x))=0 .
$$
Therefore, 
$$
\lim_{k \to +\infty}  \frac{1}{a_{k}+1}\int_0^{a_{k}+1}H(\phi(t,x)) dt =
\lim_{k\to +\infty} \frac{2c}{a_{k}+1}\sum_{j=0}^{a_k}   \chi_A(\widehat{F}^{j}(x)) +o(\varepsilon)=
2c L + o(\varepsilon)
$$
where the last equality follows from (\ref{Hist_Beha1}). Similarly,
$$
\lim_{k \to +\infty}  \frac{1}{b_{k}+1}\int_0^{b_{k}+1}H(\phi(t,x)) dt =2c \widehat{L} + o(\varepsilon). 
$$
Since $L\neq \widehat{L}$,  then for sufficiently small $\varepsilon$ and  for all $x$ in the open set $\mathcal{S}(\mathcal{U})$ we have:
$$
\lim_{a_i \rightarrow +\infty} \frac{1}{a_i} \int_0^{a_i} H(\phi(t,x)) dt   \neq \lim_{b_i \rightarrow +\infty} \frac{1}{b_i} \int_0^{b_i} H(\phi(t,x)) dt  $$

It follows that  for each $\lambda$ in $(\xi_i,\zeta_i)$ there are vector fields $g$ arbitrarily close to $f_\lambda$ in the $C^2$-topology such that there is an open set of initial conditions  for which the solution of $\dot x=g(x)$ has historic behaviour, as claimed.
 \Qed

In the proof of Theorem~\ref{teoremaHistoric}, the conditions defining the set $\mathcal C$ are only used to obtain Theorem~\ref{teorema tangency}.
Hence,  for Theorem~\ref{teoremaHistoric}, condition (P\ref{P3}) may be replaced  in the definition of $\mathcal C$ by the assumption that  two successive nodes have different chirality,
as remarked after the proof of  Theorem~\ref{teorema tangency}.

Heteroclinic tangencies also create new tangencies near them in phase space and for nearby parameter values.  
Based on \cite{Rodrigues2015, Takens94}, it should be possible to obtain a topological interpretation of the asymptotic properties of these non-converging time averages and obtain a complete set of moduli for the attracting cycle.

\section{An example}
\label{Example}

In this section we construct a family of vector fields in $\RR^3$ satisfying properties (P\ref{P1})--(P\ref{P5}). Thus, via  Theorem~\ref{teoremaHistoric}, we provide an explicit example where trajectories with historic behaviour have positive Lebesgue measure. Our example relies on Bowen's example described in \cite{Takens94}. 
This is a vector field in the plane with structurally unstable connections.
We use  the  techniques developed by Aguiar \emph{et al} \cite{ACL06, Rodrigues} combined with symmetry breaking, to lift Bowen's example to a vector field in $\RR^3$ with   periodic solutions having robust connections arising from  transverse intersections of invariant manifolds.

\subsection{The starting point}

\begin{figure}
\begin{center}
\includegraphics[height=8cm]{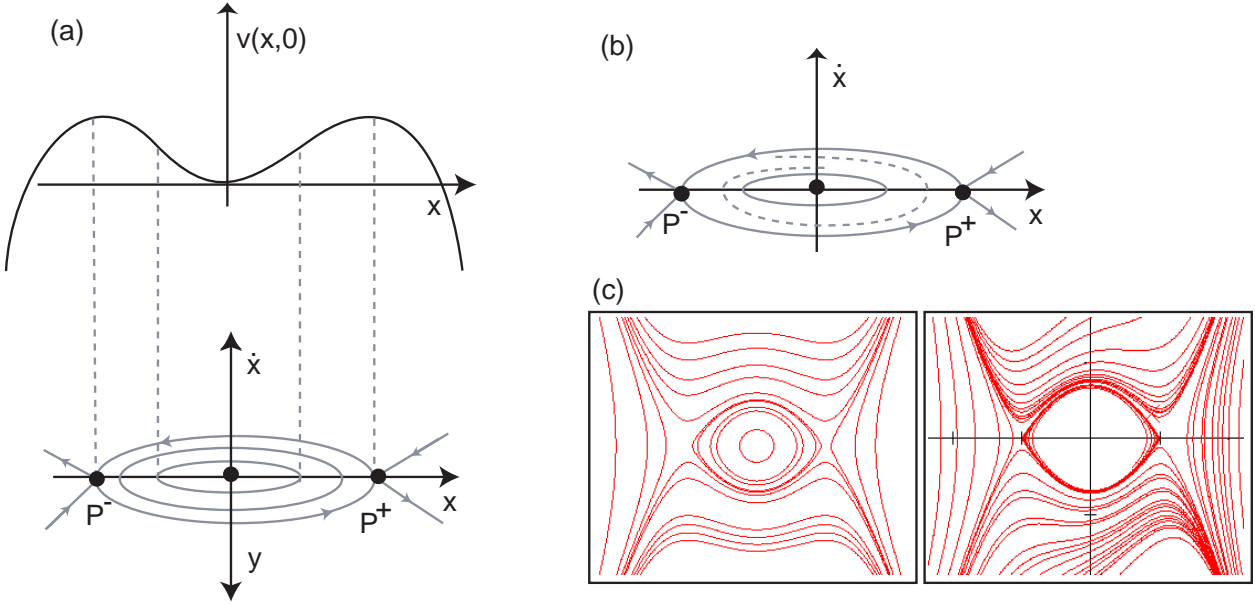}
\end{center}
\caption{\small (a): First integral and energy level of $\vv(x,0)$. (b) First perturbation. (c) Numerics for $\varepsilon=0$ and $\varepsilon =0.05$. }\label{example1}
\end{figure}

Consider the differential equation $(\dot x,\dot y)=g(x,y)$ given by
\begin{equation}
\label{example 1}
\left\{ 
\begin{array}{l}
\dot{x}=-y \\ 
\dot{y}=x-x^3
\end{array}
\right.
\end{equation}
that is equivalent to the second order equation $\ddot{x}=x-x^3$.
Its equilibria are $O = (0,0)$ and $P^\pm=(\pm 1, 0)$. 
This is a conservative system, with first integral 
$\dpt\vv(x,y)= \frac{x^2 }{2}\left(1-\frac{x^2}{2} \right)+\frac{y^2}{2}$.
From the graph of $\vv$ (see Figure~\ref{example1} (a)) it follows that the origin $O$ is a centre and the equilibria $P^\pm$ are saddles.
The equilibria $P^\pm$ are contained in the $\vv$-energy level $\vv(x,y)=1/4$ hence  there are two one-dimensional connections, one from $P^+$ to $P^-$ and another from $P^-$ to $P^+$. Denote this cycle by $\Gamma_1$. 
The region bounded by this cycle, that is  filled by closed trajectories,  will be called the \emph{invariant fundamental domain}. 
For $(x,y)\ne(0,0)$ inside the fundamental domain we have $0\le\vv(x,y)<1/4$
and the boundary of the fundamental domain intersects the $x=0$ axis at the points $(0,\pm\sqrt{2}/2)$.

\subsection{An expression for Bowen's example}\label{subsecBowen}
For a given $\varepsilon$, such that  $0<\varepsilon<\!\!<1$, consider the following  perturbation of  (\ref{example 1}):
\begin{equation}
\label{example 2}
\left\{ 
\begin{array}{l}
\dot{x}=-y \\ 
\dot{y}=x-x^3 -\varepsilon y\left(\vv(x,y)-\frac{1}{4}\right)
\end{array}
\right.
\end{equation}
 
\begin{lemma}
In the flow of equation (\ref{example 2}), the cycle $\Gamma_1$ persists and is asymptotically stable with respect to the invariant fundamental domain. 
\end{lemma}

\textbf{Proof:}
The  term $-(\vv(x,y)-1/4)$ is zero on $\Gamma_1$ and  positive in the interior of the fundamental domain. Therefore, the perturbing term $-y(\vv(x,y)-1/4)$ has the same sign as $y$.
Hence the heteroclinic connections $[P^+ \rightarrow P^-]$ and $[P^- \rightarrow P^+]$ are preserved and solutions starting away from the origin inside the fundamental domain approach the cycle when time goes to infinity as in Figure~\ref{example1} (b) and (c).
\Qed

\subsection{Translating the cycle}
For  $z^2=y+1$, Bowen's example \eqref{example 2} takes the form
\begin{equation}
\label{example 3}
\left\{
\begin{array}{l}
\dot{x}= 2z^2(1-z^2)\\ 
\dot{z}= z\left(x-x^3-\varepsilon \left(\frac{x^2}{2}- \frac{x^4}{4} + \frac{(z^2-1)^2}{2}-\frac{1}{4}\right)(z^2-1)\right)
\end{array}
\right.
\end{equation}

\begin{lemma}
The  following assertions hold for equation \eqref{example 3}:
\begin{enumerate}
\item\label{L1}
 it is  $\ZZ_2$-equivariant under the reflection on the  $z=0$ axis;
\item\label{L2}
 the $z=0$ axis is flow-invariant;
\item\label{L3}
the dynamics of \eqref{example 3}  in the $z>0$ half-plane
is orbitally equivalent to that of  \eqref{example 2} in the $y>-1$ half-plane.
\end{enumerate}
\end{lemma}

\textbf{Proof:} 
Assertion {\sl 1.} is a simple calculation, and it implies assertion {\sl 2.}
 For {\sl 3.} with $z> 0$, use $z^2=y+1$ and $\dot{z}= \frac{\dot{y}}{2z}$ to put equation \eqref{example 2} in the form: 
$$
\left\{
\begin{array}{l}
\dot{x}= 1-z^2\\ 
\dot{z}= \frac{1}{2z}\left(x-x^3-\varepsilon\left(\frac{x^2}{2}- \frac{x^4}{4} + \frac{(z^2-1)^2}{2}-\frac{1}{4}\right)(z^2-1)\right) .
\end{array} 
\right.
$$
Multiplying both equations by the positive term $2z^2$ does not affect the phase portrait and thus \eqref{example 2} with  $y>-1$ is orbitally equivalent to \eqref{example 3} with  $z>0$. 
\Qed

\subsection{The lifting}\label{secLifting}
Now we are going to use a technique presented in  \cite{ACL06, Melbourne, Rodrigues} which consists essentially in three steps:
\begin{enumerate}
\item 
Start with a vector field on $\RR^2$ with a heteroclinic cycle where $\dim Fix (\gamma)=1$, $\gamma \in \mathbf{O}(2)$.
The heteroclinic cycle involves two equilibria in $Fix (\gamma)$ and one-dimensional heteroclinic connections 
that do not intersect the line  $Fix (\gamma)$.
\item 
Lift this to a vector field on $\RR^3$ by rotating it around $Fix (\gamma)$.
This transforms  one-dimensional heteroclinic connections into two-dimensional heteroclinic connections. The resulting vector field is $\So$-equivariant  under a 3-dimensional representation of $\So$.
The attracting character of the cycle is preserved by the lifting.
\item 
Perturb the vector field to destroy the $\So$-equivariance and so that the two-dimensional heteroclinic connections perturb to transverse connections. 
\end{enumerate}

Take $(x,z,\theta)$ to be cylindrical coordinates in $\RR^3$ with radial component $z$ and let 
$(x,z_1,z_2)=(x,z\cos \theta, z \sin \theta)$ be the corresponding Cartesian coordinates.
Adding  $\dot{\theta}=1$ to  \eqref{example 3} we obtain:
\begin{equation}
\label{example 5}
\left\{
\begin{array}{l}
\dot{x}= 2(1-z_1^2-z_2^2)(z_1^2+z_2^2)\\ 
\dot{z}_1= z_1 \left[x-x^3-\varepsilon(z_1^2+z_2^2-1)\left(\frac{x^2}{2}- \frac{x^4}{4} + \frac{(z_1^2+z_2^2-1)}{2}-\frac{1}{4}\right)\right] -z_2\\
\\
\dot{z}_2= z_2 \left[x-x^3-\varepsilon(z_1^2+z_2^2-1)\left(\frac{x^2}{2}- \frac{x^4}{4} + \frac{(z_1^2+z_2^2-1)}{2}-\frac{1}{4}\right)\right]+z_1 .
\end{array}
\right.
\end{equation}

\begin{lemma}\label{lemaGamma0}
The flow of \eqref{example 5}  for $\varepsilon>0$ has a 
heteroclinic cycle $\Gamma_0$ that satisfies (P\ref{P1})--(P\ref{P3}),  consisting of  two hyperbolic closed trajectories and two surfaces  homeomorphic to  cylinders.
 The cycle  $\Gamma_0$ is asymptotically stable with respect to the lifting of the fundamental domain of 
 \eqref{example 2}.
\end{lemma}

\textbf{Proof:}
We follow the arguments of \cite{ACL06, Rodrigues}. The periodic solutions are defined by:
$$
\mathcal{P}_1: 	\quad x=1, \quad  z_1^2+z_2^2=1
\qquad
\mbox{and}
\qquad
\mathcal{P}_2: 	\quad x=-1, \quad 
z_1^2+z_2^2=1.
$$
The connections are the lift of the one-dimensional connections, rotated around  the fixed-point subspace of the symmetry.
 It follows that the heteroclinic connections are two-dimensional manifolds diffeomorphic to cylinders  and a branch of the stable manifold of each periodic solution coincides with a branch of the unstable manifold of the other.
 As remarked above, the stability of the cycle is preserved.
\Qed

\subsection{Time averages}
Theorem \ref{Main1} applied to (\ref{example 5}) says that if $\phi(t,X)\subset \mathcal{B}(\Gamma_0)$ is a non-trivial solution of the differential equation, then  the accumulation points of the time average $\frac{1}{T} \int_0^T \phi(t,X) dt $  lie in the boundary of the segment joining the points 
\begin{equation}
A_1= \left(\frac{e_2-{c_1}}{e_2+{c_1}}, 0, 0\right) \qquad \text{and} \qquad A_2= \left(\frac{e_1-{c_2}}{e_1+{c_2}}, 0, 0 \right) .
\end{equation}
Since $e_1=c_1=e_2=c_2= \sqrt{2}$, the points $A_1$ and $A_2$ coincide,  although the centres of gravity of $\mathcal{C}_1$ and $\mathcal{C}_2$ do not.
The polygon ensured by Theorem  \ref{Main1} degenerates into a single point, the origin.
This is in contrast to the example constructed in \cite{Rodrigues}where the polygon is degenerate because the centres of gravity of the nodes coincide.
 Usually, for initial conditions in the basin of attraction of  heteroclinic cycles, the time averages do not converge as $t \rightarrow +\infty$.
  However, if the vector field has symmetry, some non-generic properties appear. 
  To destroy this degeneracy  and obtain historic behaviour, it is enough to replace in \eqref{example 1} the first integral by:
 $$
 \tilde{\vv}(x,y) = -(x-1)^2(x+1)^2\left(1+\frac{x^2}{2}+x^2 \right)+\frac{y^2}{2}. 
 $$
For this case, the contracting and expanding eigenvalues at the two equilibria satisfy the conditions:
$$
\mu_1= \frac{\sqrt{5}}{\sqrt{3}}\qquad \text{and} \qquad \mu_2= \frac{\sqrt{3}}{\sqrt{5}}
$$
and thus, for the lift of the corresponding system, $A_1\neq A_2$. In particular, the Birkhoff time averages do not converge and thus they have historic behaviour. 
The next step, the second perturbation, will be performed for \eqref{example 5}, constructed using the first integral $\vv$, but it could also be done starting with $\tilde{\vv}$.

\subsection{The second perturbation}
We perturb \eqref{example 3} by adding to the equation for $\dot z_1$ a term depending on $\lambda$, as follows:
\begin{equation}
\label{example 4}
\left\{
\begin{array}{l}
\dot{x}= 2(1-z_1^2-z_2^2)(z_1^2+z_2^2)\\ 
\dot{z}_1= z_1 \left[x-x^3-\varepsilon(z_1^2+z_2^2-1)\left(\frac{x^2}{2}- \frac{x^4}{4} + \frac{(z_1^2+z_2^2-1)^2}{2}-\frac{1}{4}\right)\right] -z_2+\lambda(x^2-1)\\
\\
\dot{z}_2= z_2 \left[x-x^3-\varepsilon(z_1^2+z_2^2-1)\left(\frac{x^2}{2}- \frac{x^4}{4} + \frac{(z_1^2+z_2^2-1)^2}{2}-\frac{1}{4}\right)\right]+z_1 .
\end{array}
\right.
\end{equation}
A geometric argument is used to show that the invariant manifolds of the periodic solutions of  \eqref{example 4} intersect transversely.

\begin{lemma}
For small $\lambda>0$ and $\varepsilon>0$, 
the flow of (\ref{example 4}) has a heteroclinic cycle associated to two hyperbolic periodic solutions, $\mathcal{P}_1$ and $\mathcal{P}_2$, satisfying properties (P\ref{P1})--(P\ref{P5}). 
\end{lemma}

\begin{figure}[ht]
\begin{center}
\includegraphics[height=5cm]{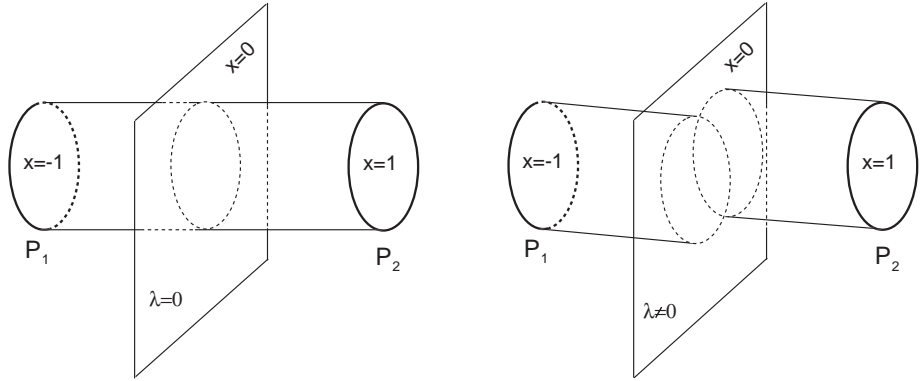}
\end{center}
\caption{\small Sketch of the invariant manifolds of $\mathcal{P}_1$ and $\mathcal{P}_2$ in  \eqref{example 4}. 
For $\lambda=0$  (left) pairs of branches of invariant manifolds of the closed trajectories coincide.
For $\lambda \neq 0$ (right) each pair of invariant manifolds meets transversely at two curves corresponding to  two points on the plane $x=0$. }
\label{intersection_example}
\end{figure}

\textbf{Proof:}
Properties (P\ref{P1})--(P\ref{P3}) follow from the construction and from Lemma~\ref{lemaGamma0}.
The perturbing term $\lambda(x^2-1)$ is zero on the planes $x=\pm 1$ that contain the cycles $\mathcal{C}_1$ and $\mathcal{C}_2$, so the periodic solutions persist.
Since $\lambda(x^2-1)$ is positive for $-1<x<1$, then when  $\lambda$ increases from zero, $W^u(\mathcal{P}_a)$, $a\in\ZZ_2$, moves towards larger values of $z_1$, while $W^s(\mathcal{P}_{a+1})$ moves in the opposite direction.
In particular, on the plane $x=0$, for $\lambda \neq 0$, each pair of invariant manifolds meets transversely at two points (Figure~\ref{intersection_example}).
Hence, there are two curves where each pair of invariant manifolds of the periodic solutions meets transversely and properties (P\ref{P4})--(P\ref{P5}) hold. 
\Qed

Let $f_\lambda$ be the family of vector fields of (\ref{example 4}). 
Theorem~\ref{teorema tangency} says that for each $a\in\ZZ_2$ there is a sequence of values of $\lambda$ for which $W^u(\mathcal{P}_a)$ is tangent to $W^s(\mathcal{P}_{a+1})$, and that for other values of $\lambda$ arbitrarily close to the connections there are closed trajectories with homoclinic tangencies.
It follows that for these values of $\lambda$ the vector field $f_\lambda$ lies in the closure of a Newhouse domain.
Theorem~\ref{teoremaHistoric} ensures that there are sequences $0<\xi_i<\zeta_i<\xi_{i+1}$, with $\lim\zeta_i=0$,
such that for each $\lambda$ in $(\xi_i,\zeta_i)$ there are vector fields $g$ arbitrarily close to $f_\lambda$ in the $C^2$-topology such that there is an open set of initial conditions  for which the solution of $\dot x=g(x)$ has historic behaviour.

    \appendix
\section{Appendix: $C^2$-Linearizing the hyperbolic periodic solution}\label{appendix}

For $a \in \{1, \ldots, k\}$, let $\Pi_a$ be a cross section transverse to the flow at $p_a \in \mathcal{P}_{a}$. 
Since $ \mathcal{P}_{a}$ is hyperbolic, there is a neighbourhood of $p_a$ where the first return map to $p_a$, denoted by $\pi_a$, is {$C^1$-conjugate} to its linear part.  Moreover: 

\begin{lemma}
Let $\pi_a$ be  the first return map to $\Pi_a$.
For each $r\geq 2$ there is an open and dense subset of $\RR^2$ such that, if the eigenvalues $(c_a,e_a)$ of  $d\pi_a$ lie in this set, then there is a neighbourhood $V^*_a$ of $p_a$ in $\Pi_a$ where $\pi_a$ is $C^r$ conjugate to its linear part.\end{lemma}

\textbf{Proof:}
Let $r\geq 2$. In order to ensure the existence of a $C^r$ conjugacy between $d\pi_a$ and the first return map to $\Pi_a$, we use the Takens' criterion \cite[Sections 1 and 5]{Takens71} which asks for the Sternberg $\alpha \left(d\pi_a, k\right)$-condition. 
Following  Takens' terminology \cite{Takens71}, let us define:
$$
\lambda_1=\lambda_c = e^{-c_a}<1, \quad \lambda_2=\lambda_e = e^{e_a}>1, \quad  s=u=1, h=2
$$
and
 $$
 \overline{M}= \lambda_e=\overline{m}>1 ,\quad \overline{N}= \lambda_c^{-1}=\overline{n}>1.
$$
In order to apply the criterion, we should define the function $\alpha \left(d\pi_a, k\right)$.
The definition will depends on an auxiliary function $\beta \left(d\pi_a, k\right)$. 
This proof will be divided in three steps:
characterisation of $\beta$, characterisation of $\alpha$ and application of the criterion.

\begin{enumerate}
\item  \textbf{The function $\beta$:} The value of $\beta \left(d\pi_a, k\right)$ is  that of the smallest $j \in \NN$ for which:
$$
\forall r<k, \qquad \overline{N} \overline{M}^r \overline{n}^{r-j}<1.
$$
In other words, $\beta \left(d\pi_a, k\right)$ is the smallest $j \in \NN$ for which:
\begin{equation}
\label{beta1}
\forall r<k, \qquad \Phi(\lambda_c, \lambda_e,r)\lambda_c^j<1, \qquad \text{where}\qquad \Phi(\lambda_c, \lambda_e,r)= \left(\frac{1}{\lambda_c}\right)^{1+r}\lambda_e^r.
\end{equation}
Thus $\beta$ depends on $d\pi_a$ through the latter's eigenvalues.
In particular, $\beta \left(d\pi_a, k\right) >1+r$ for $r<k$.  Moreover, for $j-(r+1)\in\NN$ large enough, the map $\Phi(\lambda_c, \lambda_e,r)$ increases with $r$. 
Therefore, it is sufficient to check  condition \eqref{beta1} for $r=k$.
Indeed,  the value of $\beta \left(d\pi_a, k\right)$ is that of the smallest $j \in \NN$ for which:
$$
(\lambda_c^{-1} \lambda_e)^k \lambda_c^{-1} \lambda_c^j<1\quad  \Longleftrightarrow \quad (\lambda_c^{-1})^{k+1} \lambda_c^{k+1}\lambda_e^k\lambda_c^{j-(k+1)}<1 \quad  \Longleftrightarrow \quad \lambda_e^k\lambda_c^{j-(k+1)}<1.
$$
Define $ j-(k+1)=l$. Then it is easy to see that $\beta \left(d\pi_a, k\right)$ is the smallest $j \in \NN$ such that $\lambda_e^k \lambda_c^{l}<1$. 
Taking logarithms, it follows that this is equivalent to:
$$
k\ln \lambda_e + l \ln \lambda_c<0 \quad  \Longleftrightarrow \quad l> -k \frac{\ln \lambda_e}{ \ln \lambda_c} .
$$
Since $l\in \NN$, its minimum value will be
$$
l=1 + \left[ -k \frac{\ln \lambda_e}{ \ln \lambda_c} \right]  \quad  \Longleftrightarrow \quad \beta \left(d\pi_a, k\right)=k+2 + \left[ -k \frac{\ln \lambda_e}{ \ln \lambda_c} \right],
$$
where $[x]$ represents the largest integer less than or equal to $x\in \RR$. 
\item \textbf{The function $\alpha$:}
The value of $\alpha \left(d\pi_a, k\right)$ is that of the smallest $j \in \NN$ for which:
$$
\forall r<\beta \left(d\pi_a, k\right), \qquad \overline{M} \overline{N}^r \overline{m}^{r-j}<1.
$$
In other words, $\alpha \left(d\pi_a, k\right)$ is the smallest $j \in \NN$ for which:
$$
\forall r<\beta \left(d\pi_a, k\right), \qquad \Phi(\lambda_c, \lambda_e,r)\lambda_e^{-j+1}<1, \qquad \text{with}\qquad \Phi(\lambda_c, \lambda_e,r)= \left(\frac{1}{\lambda_c}\right)^{r}\lambda_e^r.
$$
Since  $\Phi(\lambda_c, \lambda_e,r)$ increases with $r$, we would like to find the smallest $j \in \NN$ for which
$$
(\lambda_c^{-1}\lambda_e)^{\beta \left(d\pi_a, k\right)} \lambda_e^{-j+1}<1.
$$
If $j=\beta \left(d\pi_a, k\right)+1+l$, then:
$$
(\lambda_c^{-1})^{\beta \left(d\pi_a, k\right)} \lambda_e^{(\beta \left(d\pi_a, k\right)+1)}< \lambda_e^{(\beta \left(d\pi_a, k\right)+1)} \lambda_e^l
\quad  \Longleftrightarrow \quad 
 \lambda_c^{-\beta \left(d\pi_a, k\right)}<\lambda_e^l 
$$
that happens if and only if $ -\beta \left(d\pi_a, k\right)\ln (\lambda_c)<l \ln (\lambda_e)$.
Therefore, $\alpha \left(d\pi_a, k\right)= \beta \left(d\pi_a, k\right)+1+l$  where $l$ is the smallest integer $l$ such that $-\beta \left(d\pi_a, k\right)\ln (\lambda_c)<l \ln (\lambda_e)$. 
Noting that $\ln(\lambda_c)<0$ and $l \in \NN$, we have:
$$
l > -\frac{\beta \left(d\pi_a, k\right)\ln \lambda_c}{\ln \lambda_e}
\quad\mbox{with minimum value}\quad
l=1+ \left[  -\frac{\beta \left(d\pi_a, k\right)\ln \lambda_c}{\ln \lambda_e}\right]. 
$$

In conclusion, since $\ln \lambda_c=-c_a$ and $\ln \lambda_e=e_a$, 
it follows that:
$$
\beta \left(d\pi_a, k\right)=k+2+ \left[  \frac{k e_a}{c_a}\right] 
$$
and
$$
\alpha \left(d\pi_a, k\right)=\beta \left(d\pi_a, k\right)+1+l = k+4+\left[  \frac{k e_a}{c_a}\right
] + \left[ \left( k+2+\left[  \frac{k e_a}{c_a}\right]
 \right) \frac{c_a}{e_a}\right].
$$

\item \textbf{Applying the Sternberg 
condition:} 
In order to have $C^r$ conjugacy between $\pi_a$ and its linear part, the eigenvalues of $d\pi_a$ must satisfy the $\alpha \left(d\pi_a, r\right)$-condition, that we proceed to explain in this context.
For all $ \nu_1, \nu_2 \geq 0 $ such that $ 2\leq  \nu_1+\nu_2 \leq \alpha \left(d\pi_a, r\right)$ we should have:
$$ \lambda_c^{\nu_1-1}\lambda_e^{\nu_2} \neq 1, \qquad \lambda_e^{\nu_2-1}\lambda_c^{\nu_1}\neq 1\quad \text{and}\quad |\lambda_c^{\nu_1} \lambda_e^{\nu_2}|\neq 1
$$
Indeed, $ \lambda_c^{\nu_1} \lambda_e^{\nu_2}= e^{-\nu_1 c_a}e^{-\nu_2 e_a} =1$ if and only if $-\nu_1 c_a=\nu_2 e_a$. In summary, for all $ \nu_1, \nu_2 \geq 0 $ such that $ 2\leq  \nu_1+\nu_2 \leq \alpha \left(d\pi_a, r\right)$, the following conditions should hold: 
\begin{itemize}
\item $(\nu_1-1) c_a \neq \nu_2 e_a$
\item $(\nu_1) c_a \neq (\nu_2-1) e_a$
\item $\nu_1 c_a \neq \nu_2 e_a$.\Qed
\end{itemize}
\end{enumerate}

The set of smooth vector fields that satisfy the Sternberg $\alpha \left(d\pi_a, r\right)$-condition, for each $r\ge 2$, is open and dense in the set of vector fields satisfying  (P\ref{P1}) -- (P\ref{P5}). Hence, generically the assumptions are satisfied.

\section{Control of flight times}
\label{appendixB}

\subsection{Proof of Lemma \ref{Equalities}}
\label{appendixB1}

\begin{enumerate}
\item If $n=0$ ($n$ corresponds to the number of loops around the cycle $\Gamma_0$), it is trivial. 
For $n \geq 1$, we may write the following equality, omitting the dependence on $X$:
\begin{eqnarray*}
T_{a + nk}=&T_a &+\quad \tau_a +\tau_{a+1}+ \ldots \tau_{a+k-1}+\\
&&+\quad \tau_{a+k} +\tau_{a+k+1}+ \ldots \tau_{a+2k-1}+\ldots\\
&&+\quad \tau_{a+(n-1)k} +\tau_{a+(n-1)k+1}+ \ldots \tau_{a+nk-1}\\
\end{eqnarray*}
Using Corollary \ref{Lemma3}, the previous equality yields:
\begin{eqnarray*}
T_{a + nk}=&T_a &+ \mu_a \tau_{a-1} +\mu_a\mu_{a+1}\tau_{a-1}+ \ldots \prod_{l=0}^{k-1}\mu_{a+l}\tau_{a-1}+\\
&&+ \delta \mu_a \tau_{a-1}+\delta\mu_a\mu_{a+1}\tau_{a-1}+ \ldots \delta \left(\prod_{l=0}^{k-1}\mu_{a+l}\right)\tau_{a-1}+\ldots\\
&&+ \delta^{n-1} \mu_a \tau_{a-1}+\delta^{n-1}\mu_a\mu_{a+1}\tau_{a-1}+ \ldots \delta^{n-1}\left(\prod_{l=0}^{k-1}\mu_{a+l}\right)\tau_{a-1}=\\
=& T_a& + \frac{\delta^n-1}{\delta-1} \left(\mu_a +\mu_a \mu_{a+1} + \ldots + \prod_{l=0}^{k-1}\mu_{a+l}\right)\tau_{a-1}
\end{eqnarray*}

\item This item follows from Corollary \ref{Lemma3}. Indeed, we have:
$$
\tau_{a + nk}(X)=\mu_{a }\tau_{a+nk-1}(X)=\mu_{a}\mu_{a-1}\tau_{a+nk-2}(X)=\ldots=\delta^n\tau_{a}(X)=\delta^n\mu_a\tau_{a-1}(X).
$$
\Qed
\end{enumerate}

\subsection{Proof of Proposition \ref{Prop6}}
\label{appendixC}

We divide the proof in two lemmas. 
First we show in Lemma~\ref{lematrivial?} that it is sufficient to consider the limit when $n\to\infty$ of the averages over one turn around $\Gamma_0$.
Then in Lemma~\ref{Lemma4} we show that these averages tend to $A_a$.

\begin{lemma}\label{lematrivial?}
Let $T_\ell $, $\ell\in\NN$, be a sequence $0=T_0<T_\ell <T_{\ell +1}$ with $\lim_{\ell \to \infty} T_\ell =\infty$.
Given $h:\RR\rightarrow\RR^m$ an integrable map,
$$
\mbox{if }
\lim_{\ell \to\infty}\frac{1}{T_{\ell +1}-T_{\ell }} \int_{T_{\ell }}^{T_{\ell +1}} h(t)dt=\omega,
\qquad\mbox{then}\qquad
\lim_{\ell \to\infty}\frac{1}{T_{\ell }} \int_{0}^{T_{\ell }} h(t)dt=\omega .
$$
\end{lemma}
\textbf{Proof:}
First note that
$$
\frac{1}{T_{\ell }} \int_{0}^{T_{\ell }} h(t)dt - \omega=\frac{1}{T_{\ell }} \int_{0}^{T_{\ell }} (h(t)-\omega) dt=
\sum_{j=1}^\ell  \left(\frac{T_j-T_{j-1}}{T_\ell }\right) \left[
\frac{1}{T_j-T_{j-1}} \int_{T_{j-1}}^{T_{j}} (h (t)-\omega)dt
\right] .
$$
From the hypothesis, given $\varepsilon>0$ there exists $N_1$ such that  $\ell >N_1$ implies
$$
 \frac{1}{T_\ell -T_{\ell -1}}\left| \int_{T_{\ell -1}}^{T_{\ell }} (h (t)-\omega)dt\right|<\frac{\varepsilon}{2}.
 $$
 Let $A=\displaystyle \left|\int_0^{T_{N_1} }(h (t)-\omega) dt\right|$. 
 Since $T_\ell \to\infty$ then there exists $N_2$ such that $T_{N_2}> 2A/\varepsilon$.
Let $N_0=\max\left\{ N_1,N_2\right\}$.
If $\ell >N_0$ then
\begin{eqnarray*}
\left|\frac{1}{T_{\ell }} \int_{0}^{T_{\ell }} (h (t)-\omega )dt\right|&\le&
\frac{1}{T_{\ell }}\left| \int_{0}^{T_{N_1}} (h (t)-\omega )dt\right| +
\frac{1}{T_{\ell }} \left|\int_{T_{N_1}}^{T_{\ell }} (h (t)-\omega )dt\right|\\
&\le&\frac{A}{T_\ell } + 
\sum_{j=N_1}^\ell  \left(\frac{T_j-T_{j-1}}{T_\ell }\right) \  
\frac{1}{T_j-T_{j-1}}\left| \int_{T_{j-1}}^{T_{j}} (h (t)-\omega )dt\right| \\
&\le&\frac{\varepsilon}{2}+\sum_{j=N_1}^\ell  \left(\frac{T_j-T_{j-1}}{T_\ell }\right) \frac{\varepsilon}{2}
\le \frac{\varepsilon}{2}\left( 1+\sum_{j=1}^\ell  \frac{T_j-T_{j-1}}{T_\ell } \right)= \varepsilon .
\end{eqnarray*}
\Qed

\begin{lemma}
\label{Lemma4}
Let $f_0$ be a vector field  in $\RR^3$  satisfying (P\ref{P1})--(P\ref{P3}).  
For each $a\in\ZZ_k$, and for each $X\in\mathcal{B}(\Gamma_0)$, 
the limit  of the spatial average of of $\phi(t,X)$  over one full turn around the heteroclinic cycle $\Gamma_0$ starting at $In(\mathcal{P}_j)$ is  $A_j$. More precisely:
$$
\lim_{n\to\infty} \frac{1}{T_{a+(n+1)k}(X)-T_{a+nk}(X)}\int_{T_{a+(n+1)k}(X)}^{T_{a+nk}(X)}\phi(t, X) dt =A_a.
$$
\end{lemma}

\textbf{Proof:}
First, recall that we are assuming that, for all $a\in\ZZ_k$, the jumps from $Out(\mathcal{P}_a)$ to $In(\mathcal{P}_{a+1})$ are instantaneous (see Remark~\ref{rkFlightTimes}).
 Since $X\in\mathcal{B}(\Gamma_0)$, then for $t\in\left[T_{a+nk},T_{a+1+nk} \right]$ with large $n$, the trajectory $\phi(t,X)$ gets very close to $\mathcal{P}_a$.
Therefore, omitting the $(X)$ for shortness, we have:
 \begin{equation}
 \label{average1}
 \lim_{n\to \infty} \frac{1}{T_{a+1+nk}-T_{a+nk}}\int_{T_{a+nk}}^{T_{a+1+nk}}\phi(t, X) dt = 
 \lim_{n\to \infty}  \frac{1}{\tau_{a+nk}} \int_{T_{a+nk}}^{T_{a+1+nk}}\phi(t, X) dt = \overline{x_a}.
 \end{equation}

Without loss of generality, from now on we take $a=1$.
Then
\begin{eqnarray*}
&&
 \frac{1} {T_{k+1+nk}-T_{1+nk}} \int_{T_{1+nk}}^{T_{k+1+nk}} \phi(t, X)dt \\
&= &
 \frac{1} {T_{k+1+nk}-T_{1+nk}}\left[\int_{T_{1+nk}}^{T_{2+nk}} \phi(t, X)dt+\int_{T_{2+nk}}^{T_{3+nk}} \phi(t, X)dt +\cdots+ \int_{T_{k+nk}}^{T_{k+1+nk}} \phi(t, X)dt\right]\\
 &=&
 \sum_{b=1}^k  \frac{T_{b+1+nk}-T_{b+nk}}{T_{k+1+nk}-T_{1+nk}} 
\left[  \frac{1}{T_{b+1+nk}-T_{b+nk}}
 \int_{T_{b+nk}}^{T_{b+1+nk}} \phi(t, X)dt \right] .
\end{eqnarray*}
Recall that $T_{2+nk}-T_{1+nk}=\tau_{1+nk}$, and
by Corollary~\ref{Lemma3}, for any $b\in\{2,\ldots,k\}$
\begin{eqnarray*}
 \frac{T_{b+1+nk}-T_{b+nk}} {T_{k+1+nk}-T_{1+nk}}&=&
 \frac{\tau_{b+nk} }{\tau_{1+nk}+\tau_{2+nk}+\ldots+ \tau_{k+nk}}\\
&=&\frac{\mu_b \mu_{b-1}\ldots \mu_2 \tau_{1+nk}}{\tau_{1+nk}+\mu_2 \tau_{1+nk} +\ldots+ \mu_k \mu_{k-1}\ldots \mu_2 \tau_{1+nk}}\\
&=&\frac{\mu_b \mu_{b-1}\ldots \mu_2 }{den A_1} .
\end{eqnarray*}
Therefore, the value of
$$
\lim_{n \rightarrow +\infty}
 \frac{1} {T_{k+1+nk}-T_{1+nk}} \int_{T_{1+nk}}^{T_{k+1+nk}} \phi(t, X)dt 
$$
is, by \eqref{average1},
\begin{eqnarray*}
&&
\lim_{n \rightarrow +\infty}
\frac{1}{den A_1}\left[ \frac{1}{T_{2+nk}-T_{1+nk}}\int_{T_{1+nk}}^{T_{2+nk}} \phi(t, X)dt+
  \sum_{b=2}^k \frac{\mu_b \mu_{b-1}\ldots \mu_2}{T_{b+1+nk}-T_{b+nk}}
 \int_{T_{b+nk}}^{T_{b+1+nk}} \phi(t, X)dt
 \right]\\
 &=& \frac{ \overline{x}_{1} + \mu_2  \overline{x}_{2} + \ldots+ \mu_k \mu_{k-1}\ldots \mu_2  \overline{x}_{k}}{1+\mu_2  +\ldots+ \mu_k \mu_{k-1}\ldots \mu_2}=A_1 .
\end{eqnarray*}
\Qed

\subsection{Proof of Lemma \ref{Colinear}}

Expanding $den(A_{a})$ and $den(A_{a+1})$, yields:
$$
den(A_a) = {1+\mu_{a+1} +\mu_{a+1}\mu_{a+2}+\cdots+  \prod_{l=1}^{k-1}\mu_{l+a}}
\qquad
den(A_{a+1}) = {1+\mu_{a+2} +\mu_{a+2}\mu_{a+2}+\cdots+  \prod_{l=1}^{k-1}\mu_{l+a+1}}
$$
hence, since $ \prod_{l=0}^{k-1}\mu_{l+a+1}=\delta$, 
\begin{eqnarray*}
\mu_{a+1} den(A_{a+1}) &=& {\mu_{a+1}+\mu_{a+1}\mu_{a+2} +\mu_{a+1}\mu_{a+2}\mu_{a+2}+\cdots+  \mu_{a+1} \prod_{l=1}^{k-1}\mu_{l+a+1}}\\
&=&  \mu_{a+1}+\mu_{a+1}\mu_{a+2} +\mu_{a+1}\mu_{a+2}\mu_{a+2}+\cdots+ \prod_{l=1}^{k-1}\mu_{l+a} + \delta\\
&=& den (A_a)-(1-\delta) .
\end{eqnarray*}
For $num(A_{a})$ and $num(A_{a+1})$ we obtain
$$
num(A_a) = \overline{x}_{a} + \mu_{a+1} \overline{x}_{a+1} + \mu_{a+1}\mu_{a+2}  \overline{x}_{a+2} + \cdots+ \left(\prod_{l=1}^{k-1}\mu_{l+a} \right) \overline{x}_{a+k-1} 
$$
$$
num(A_{a+1}) = \overline{x}_{a+1} + \mu_{a+2} \overline{x}_{a+2} + \mu_{a+2}\mu_{a+3}  \overline{x}_{a+3} + \cdots+ \left(\prod_{l=1}^{k-1}\mu_{l+a+1} \right) \overline{x}_{a+k} 
$$
and, since $\overline{x}_{a+k}=\overline{x}_a$, we get
\begin{eqnarray*}
\mu_{a+1} num(A_{a+1}) &=& \mu_{a+1} \overline{x}_{a+1} + \mu_{a+1} \mu_{a+2} \overline{x}_{a+2} + \mu_{a+1} \mu_{a+2}\mu_{a+3}  \overline{x}_{a+3} + \cdots+ \mu_{a+1} \left(\prod_{l=1}^{k-1}\mu_{l+a+1} \right) \overline{x}_{a+k}\\
&=& \mu_{a+1} \overline{x}_{a+1} + \mu_{a+1} \mu_{a+2} \overline{x}_{a+2} + \mu_{a+1} \mu_{a+2}\mu_{a+3}  \overline{x}_{a+3} + \cdots +  \left(\prod_{l=1}^{k-1}\mu_{l+a} \right) \overline{x}_{a+k-1}+ \delta \overline{x}_{a}\\
&=& num(A_a) - (1-\delta)\overline{x}_a 
\end{eqnarray*}
and the lemma is proved.
\Qed

\end{document}